\pdfoutput=1
\documentclass[11pt]{amsart}
\usepackage[utf8]{inputenc}
\usepackage{amssymb,amsmath,amsthm,enumerate,enumitem,colonequals,mlmodern,tikz-cd,microtype,stmaryrd}
\usepackage[cal=euler,bfcal,bb=px,bfbb]{mathalpha}
\usepackage[top=3.75cm, bottom=3cm, left=3.5cm, right=3.5cm]{geometry}
\makeatletter
\@namedef{subjclassname@2020}{\textup{2020} Mathematics Subject Classification}
\makeatother
\usepackage{mathtools}
\usepackage{needspace}
\usepackage{xcolor}
\colorlet{darkblue}{blue!55!black}
\colorlet{darkcyan}{cyan!50!black}
\colorlet{darkgreen}{green!60!black}

\usepackage{hyperref}
\hypersetup{
    colorlinks=true,
    linkcolor= darkblue,
    urlcolor= darkcyan,
    citecolor= darkgreen,
}

\def\eqref#1{\textcolor{darkblue}{(\ref{#1})}}

\usepackage{mathrsfs}
\usepackage[pagewise]{lineno}
\let\oldequation\equation
\let\oldendequation\endequation
\renewenvironment{equation}{\linenomathNonumbers\oldequation}{\oldendequation\endlinenomath}
\expandafter\let\expandafter\oldequationstar\csname equation*\endcsname
\expandafter\let\expandafter\oldendequationstar\csname endequation*\endcsname
\renewenvironment{equation*}{\linenomathNonumbers\oldequationstar}{\oldendequationstar\endlinenomath}
\let\oldalign\align
\let\oldendalign\endalign
\renewenvironment{align}{\linenomathNonumbers\oldalign}{\oldendalign\endlinenomath}
\expandafter\let\expandafter\oldalignstar\csname align*\endcsname
\expandafter\let\expandafter\oldendalignstar\csname endalign*\endcsname
\renewenvironment{align*}{\linenomathNonumbers\oldalignstar}{\oldendalignstar\endlinenomath}


\newcounter{intro}
\newcounter{HypCounter}

\newtheorem{introthm}[intro]{Theorem}

\theoremstyle{plain}
\newtheorem{theorem}{Theorem}[section]
\newtheorem{lemma}[theorem]{Lemma}
\newtheorem{corollary}[theorem]{Corollary}
\newtheorem{proposition}[theorem]{Proposition}

\theoremstyle{definition}
\newtheorem{conjecture}[theorem]{Conjecture}

\newtheorem{example}[theorem]{Example}

\newtheorem*{hypothesis*}{Hypothesis}

\newtheorem{remark}[theorem]{Remark}

\newtheorem*{ack}{Acknowledgements}

\numberwithin{equation}{section}
\numberwithin{theorem}{section}

\makeatletter
\renewcommand{\@settitle}{\begin{center}
  \fontencoding{T1}\fontfamily{ppl}\fontseries{m}\fontshape{sc}\fontsize{18}{23}\selectfont
  \@title\par
\end{center}}
\makeatother

\title[Half-Bounded Tensor Triangular Geometry]{Spectra, Growth, and Localization in Half-Bounded Tensor Triangular Geometry}

\author[Liangyu Zhu]{Liangyu Zhu}
\address{School of Mathematical Sciences, Anhui University, Hefei 230601, China}
\email{a126112006@stu.ahu.edu.cn}

\date{October 25, 2026}



\begin{document}
\begin{abstract}
Over a noetherian scheme, classical cohomological support determines
thick tensor ideals of perfect complexes. In half-bounded derived categories, complexes with
the same support can generate different radical thick tensor ideals.
We study how these distinctions affect the Balmer spectrum and
localization. For a commutative noetherian ring $R$, we prove that
$\operatorname{Spc}(\mathbf D^-_{\mathrm{fg}}(R))$ is noetherian if and
only if $R$ is artinian. For a nonempty noetherian scheme $X$, the
corresponding condition is $\dim X=0$. For such a ring $R$ and a finite
group $G$, $\operatorname{Spc}(\mathbf D^-_{\mathrm{fg}}(RG))$ is
noetherian if and only if $R$ is artinian and $|G|$ is invertible in
$R$. Using the Sanders--Zhang classification over a discrete valuation
ring, we compare objects whose annihilator exponents grow at different
rates in different degrees. We lift the relevant generation relations
before changing coefficients and obtain spectral retractions. These
comparisons persist after thick quotients and algebraic localization
under explicit uniform annihilation and nonvanishing conditions. A
second argument uses nonnilpotent classes of positive degree and passes
from cyclic subgroups to finite groups. Applications include relative
results for schemes and bounded-below spectra with degreewise finitely
generated homotopy. We also describe the cohomological information shared
by the resulting ideals. The localization comparisons for bounded-below
complexes of finitely generated modules and for bounded-below derived
objects with finitely generated cohomology are established separately.
\end{abstract}
\maketitle
\newpage
\setcounter{tocdepth}{1}\tableofcontents\setcounter{tocdepth}{2}
\clearpage
\section{Introduction}\label{sec:intro}

Localization and support are basic tools in tensor triangular geometry.
They turn questions about generating objects into questions about a
topological space: where is an object supported, which objects can be
built from it, and what remains after some of them have been removed?
The classical examples come from algebraic geometry, stable homotopy
theory and modular representation theory
\cite{balmer2005spectrum,balmer2010spectra}.

For a noetherian scheme $X$, the theorems of Hopkins, Neeman and
Thomason classify thick tensor ideals of
$\mathbf D^{\mathrm{perf}}(X)$ by specialization-closed subsets of $X$
\cite{hopkins1987global,neeman1992chromatic,Thomason}. If $P$ and $Q$
are perfect complexes, then
\[
 Q\in\operatorname{Thick}^{\otimes}(P)
 \quad\Longleftrightarrow\quad
 \operatorname{Supp}Q\subseteq\operatorname{Supp}P.
\]
Balmer recovers the underlying space from prime tensor ideals
\cite[Corollary~5.6]{balmer2005spectrum}:
\[
 \operatorname{Spc}(\mathbf D^{\mathrm{perf}}(X))\simeq X.
\]
Thus classical cohomological support answers the generation question
for perfect complexes. In this setting a noetherian space remains
noetherian when recovered from its tensor category.

The distinction between \emph{rigid tensor triangular geometry}
and \emph{non-rigid tensor triangular geometry} is central to our
questions. In a rigid tensor triangulated category, every object $M$
has a tensor dual $M^\vee$. The evaluation and coevaluation identities
make $M$ a retract of $M\otimes M^\vee\otimes M$. Iterating this
observation gives
\[
 M^{\otimes n}\in\mathcal J\quad\Longrightarrow\quad M\in\mathcal J
 \qquad(n\ge1)
\]
for every thick tensor ideal $\mathcal J$. Thus every thick tensor
ideal is radical
\cite[Remark~4.3 and Proposition~4.4]{balmer2005spectrum}.
The categories of perfect complexes and of finite spectra are rigid,
as is the bounded derived category of finite-dimensional
representations over a field. Duals allow information about tensor
powers to be brought back to the object itself.

In non-rigid tensor triangular geometry, objects without tensor duals
are allowed, and this argument is no longer available for them. We
must distinguish a thick tensor ideal from its radical. Balmer's
classification of radical thick tensor ideals remains valid
\cite[Theorem~4.10]{balmer2005spectrum}; what changes is the category
to which it is applied, and hence the spectrum and the generation
relations that it records. Over a noetherian scheme, the dualizable
objects in $\mathbf D^-_{\mathrm{Coh}}(X)$ are precisely the perfect
complexes \cite[Tag~0FPD]{stacks-project}. Enlarging the category
therefore changes the objects among which generation is tested.
An object with an infinite cohomological tail is not dualizable.
On a singular scheme, a bounded coherent complex can already fail
to be perfect.

It would be too quick to attribute every topological difference to
the absence of duals. For $R=k[\varepsilon]/(\varepsilon^2)$, the
module $k$ has an infinite periodic free resolution and is not
perfect, but $\operatorname{Spc}(\mathbf D^-_{\mathrm{fg}}(R))$ is
a single point by the artinian case of
\cite[Theorem~6.5]{Matsui_2017}. In the other direction, finite
$p$-local spectra form a rigid category with a nonnoetherian Balmer
spectrum \cite[Section~9]{balmer2010spectra}. Neither rigidity nor
its failure settles the noetherianity question. Which generation
relations are responsible in the categories considered here? The
perfect-complex case gives us a clear starting point, because both
the geometry and the generation relations are known.

There is also an extensive theory for categories admitting arbitrary
coproducts. Tensor idempotents describe localization in the rigidly
compactly generated setting \cite{balmerfavi2011idempotents}; ring
actions and stratification give support classifications under their
respective hypotheses \cite{BIK2011}. Actions on triangulated module
categories extend these methods further \cite{Stevenson}. In the
categories considered here, we retain a different finiteness
condition: cohomology is coherent, or finitely generated, in each
degree, while only one end of the complex is bounded. These categories
contain infinite cohomological tails but need not contain arbitrary
coproducts. The subcategories we study are thick, so their generation
uses finite sums, cones and retracts, together with the given action.
This finiteness condition is important. In a category closed under
all coproducts, a localizing subcategory can contain an object
because it contains all pieces in an infinite construction.
Membership in a thick subcategory requires a finite construction.
The support theorems for localizing subcategories do not by
themselves tell us which infinite complexes can be obtained by
finitely many such operations. We keep the cohomological finiteness
that occurs naturally in algebra and geometry, and ask that
generation question inside the resulting small categories.

We would like to understand what happens when dualizability is lost
in this way. Passing from $\mathbf D^{\mathrm{perf}}(X)$ to
$\mathbf D^-_{\mathrm{Coh}}(X)$ retains a tensor triangulated
category, but an object with infinitely many nonzero cohomology
sheaves is no longer dualizable. Ordinary support still makes sense.
Every bounded truncation is still available. How much do these data
tell us about the whole object?

\emph{When does the Balmer spectrum remain noetherian after passing
to a half-bounded category? How much of a tensor ideal is determined
by its bounded truncations or by its cohomology?}

One might hope to answer the second question by assembling the
bounded truncations. The difficulty is that a thick subcategory is
not required to perform an infinite assembly. Even if every bounded
truncation can be built from a given collection, the number of
extensions or tensor factors needed may grow with the truncation.
There may be no single finite construction for the whole tail. This
is where we expect information beyond ordinary support to occur.
Our concern is to find that information explicitly and to understand
which operations preserve it.

The existing work over commutative noetherian rings already shows
that this is a real issue. Matsui--Takahashi classify the thick tensor
ideals of $\mathbf D^-_{\mathrm{fg}}(R)$ generated by bounded
complexes, equivalently by perfect complexes, and distinguish them
from the other tensor ideals \cite{Matsui_2017}. They show that the
Balmer spectrum is $\operatorname{Spec}R$ when $R$ is artinian
\cite[Theorem~6.5]{Matsui_2017}. Matsui then proves the necessary
condition \cite[Theorem~3.1]{matsui2018connectedness}
\[
 \operatorname{Spc}(\mathbf D^-_{\mathrm{fg}}(R))\text{ noetherian}
 \quad\Longrightarrow\quad \operatorname{Spec}R\text{ finite}
 .
\]
But finiteness of $\operatorname{Spec}R$ still allows a nontrivial
specialization: a discrete valuation ring has exactly two prime
ideals. Could such a small classical spectrum nevertheless carry
infinitely many different half-bounded generation relations?

Sanders--Zhang answer this question by computing the DVR spectrum
\cite[Theorems~10.14 and~11.1]{sanders2026balmer}. Let $V$ be a DVR
with uniformizer $\pi$ and maximal ideal $\mathfrak m$. For a positive
nondecreasing function $f$, consider
\[
 T_f=\bigoplus_{n\ge0}V/\pi^{f(n)}V[n],\qquad
 J_f=\sqrt{\operatorname{Thick}^{\otimes}(T_f)}.
\]
These objects all have support $\{\mathfrak m\}$. Their ideals have
the same bounded and perfect parts, by
\cite[Proposition~2.9 and Corollary~2.11]{Matsui_2017}:
\[
\begin{aligned}
 J_f\cap\mathbf D^b_{\mathrm{fg}}(V)
   &=\mathbf D^b_{\mathrm{fg}}(V)_{\{\mathfrak m\}},\\
 J_f\cap\mathbf D^{\mathrm{perf}}(V)
   &=\mathbf D^{\mathrm{perf}}(V)_{\{\mathfrak m\}}.
\end{aligned}
\]
Nevertheless, $f(n)=\lceil(n+1)^\alpha\rceil$ gives strictly
increasing radical ideals as $\alpha>0$ increases. The resulting
half-bounded spectrum is nonnoetherian. The two-point spectrum of
perfect complexes and this larger spectrum are therefore very
different.

What has ordinary support forgotten? It records where torsion
occurs, but not how rapidly the annihilator exponents increase along
the complex. Tensor generation can distinguish those rates. We take
this as a reason to study non-rigid tensor triangular geometry by
following explicit objects and the finite constructions relating
them. Balmer support continues to classify radical tensor ideals
\cite[Theorem~4.10]{balmer2005spectrum}; the point is to understand
the additional distinctions in that support space which classical
cohomological support does not see. There are two different losses
of information to keep apart. For a thick tensor ideal over $R$,
\[
 \mathcal J\subseteq\sqrt{\mathcal J}
 \subseteq
 \{X:\operatorname{Supp}X\subseteq\operatorname{Supp}\mathcal J\}.
\]
Balmer support recovers the middle term. Ordinary cohomological
support specifies the last term, the tame closure of
Matsui--Takahashi \cite[Lemma~5.8]{Matsui_2017}. Our examples already
compare radical ideals, so their differences remain after passing
to the middle term. Can we describe these differences over a
general ring, rather than only over a DVR? Can we retain them in
the other non-rigid categories that arise from geometry and
representations?

Can the DVR calculation be used over a general ring? A map
$u:R\to V$ is a natural way to try, and Sanders--Zhang discuss this
reduction in \cite[Remark~6.27]{sanders2026balmer}. There are two
difficulties. Such a map need be neither flat nor finite, and an
element $r\in R$ whose image is a power of $\pi$ may be a zero
divisor. Replacing $V/\pi^nV$ by $R/r^nR$ therefore need not give
the required derived base change. Moreover, a generation relation
after applying $-\otimes_R^{\mathbf L}V$ need not come from one
before applying it. Lifting the individual objects does not
automatically lift the relations between their generated ideals.
For the spectral embedding below it will be enough to lift
principal radical ideals. To preserve intersections and obtain a
spectral retraction, we need explicit generation relations before
changing coefficients. A functor preserves these relations in the
forward direction, but does not supply the missing constructions in
the original category.

We use two-term free complexes instead of cyclic modules. If
$u(r)=\varepsilon\pi^c$, with $\varepsilon\in V^\times$ and $c>0$,
then
\[
 \operatorname{Cone}(r^a:R\to R)\otimes_R^{\mathbf L}V
       \simeq V/\pi^{ca}V.
\]
Taking the degreewise finite sums of their shifts gives the objects
needed for the principal radical ideals in the DVR classification. Multiplication of all exponents
by $c$ leaves the principal radical ideal over $V$ unchanged, by the
Sanders--Zhang classification. Thus each principal radical ideal
can be lifted, and this already gives the spectral embedding used
in the ring noetherianity criterion.

We also want to know more than whether the DVR spectrum embeds.
Can the lifts be chosen so that intersections and generated unions
are preserved? Can arbitrary radical ideals be transported back,
and can the embedding have a spectral retraction? These questions
require generation relations in the original category. We prove
them by finite tensor, cone and retract constructions, using the
DVR methods of \cite[Sections~7 and~10]{sanders2026balmer}. The
intersection calculation is essential for the stronger conclusions.

For positive nondecreasing functions we write
\[
 f\le_\mu g\quad\Longleftrightarrow\quad
 f(n)\le A g(2^kn)\ (n\gg0)
 \quad\text{for some }A\ge1,\ k\ge0.
\]
In a tensor triangulated category containing the indicated sums, let
\[
 K_f(s)=\bigoplus_{n\ge0}
       \operatorname{Cone}(s^{f(n)}:\mathbf1\to\mathbf1)[n],
 \qquad J_f=\sqrt{\operatorname{Thick}^{\otimes}(K_f(s))}.
\]

\begin{introthm}\label{intro:valuation-retraction}
Let $\mathcal K$ be an essentially small full thick tensor subcategory
of a tensor triangulated category with countable coproducts preserved
by tensoring. Suppose $\mathcal K$ contains
$\bigoplus_{n\ge0}\mathbf1[n]^{d_n}$ for all finite $d_n$, and let
$s\in\operatorname{End}(\mathbf1)$.
Then $[f]_\mu\mapsto J_f$ preserves order, joins and meets:
\[
 f\le_\mu g\ \Longrightarrow\ J_f\subseteq J_g,
 \qquad J_f\vee J_g=J_{\max(f,g)},\qquad
 J_f\cap J_g=J_{\min(f,g)}.
\]
If an exact symmetric monoidal functor
$F:\mathcal K\to\mathbf D^-_{\mathrm{fg}}(V)$ preserves
$\bigoplus_{n\ge0}\mathbf1[n]$ and satisfies
$F(s)=\varepsilon\pi^c$, with $\varepsilon\in V^\times$ and $c>0$,
then
\[
 J_f\subseteq J_g\quad\Longleftrightarrow\quad f\le_\mu g.
\]
The induced map on radical ideals has a section preserving arbitrary
joins and finite meets, and $\operatorname{Spc}(F)$ has a spectral
retraction.
\end{introthm}

Lemma~\ref{lem:koszul-growth-bound} gives these constructions, allowing
finitely many summands in each degree. They prove generation in
$\mathcal K$ and the intersection formula; tensoring with $V$ shows
when the inclusions are strict. After the lifts have been chosen to
preserve joins and intersections, ideal completion and Hochster
duality give the section on all radical ideals and the spectral
retraction
\cite[Sections~1.2 and~3]{kockpitsch2017hochster}.
Theorem~\ref{thm:valuation-growth-retraction} gives the details.

There is a further difficulty when we localize. Over a commutative
noetherian ring $R$, let $\Phi\subseteq R$ be multiplicative and
compare
\[
 \exists\phi\in\Phi:\ \phi\operatorname{id}_M=0,
 \qquad
 \Phi^{-1}H^iM=0\quad\text{for every }i.
\]
The first asks for one scalar to kill the whole object; the second
allows a different annihilator in each degree. For bounded complexes
with finitely generated cohomology, finite generation and a finite
Postnikov filtration give their equivalence. For a half-bounded
complex, the required exponents can grow without bound. Over the
DVR $V$, for example,
\[
 E=\bigoplus_{n\ge0}V/\pi^{n+1}V[n],\qquad
 E\otimes_V\operatorname{Frac}(V)=0,\qquad
 \pi^a\operatorname{id}_E\ne0\quad(a\ge0).
\]
The last assertion can be checked on $H^{-a}E=V/\pi^{a+1}V$.
Every bounded truncation has a uniform annihilator; the whole
complex does not.

\emph{How does this failure affect localization and the Balmer
spectrum? How large exactly is the difference between geometric
localization and algebraic localization, and how can we describe
this difference?}

These are the localization questions we pursue. Algebraic central
localization inverts $\Phi$ in morphism groups
\cite[Corollary~3.10]{balmer2010spectra}. Its kernel is the
subcategory of objects killed by a single scalar in $\Phi$.
Geometric localization means extension of scalars to $\Phi^{-1}R$;
its kernel is exactly the degreewise torsion subcategory. Hence the two
localizations treat $E$ differently. This also explains why the
usual perfect-complex localization statement is an insufficient
guide here: for perfect complexes the two constructions agree
after taking direct summands, whereas in the half-bounded category
they can differ even on whether an object is zero
\cite[Remarks~6.21--6.24]{sanders2026balmer}.

The $t$-structure tells us what the algebraic quotient has retained.
Writing $\mathcal S^-_\Phi$ for uniformly annihilated objects and
$\mathcal N^-_\Phi$ for degreewise torsion objects, we have
\[
 \operatorname{PC}^-(\mathcal S^-_\Phi)=\mathcal N^-_\Phi.
\]
Here Postnikov closure means closure under truncations and passage
from all bounded standard truncations to their half-bounded object.
Proposition~\ref{prop:central-comparison} identifies the comparison
kernel with $\mathcal N^-_\Phi/\mathcal S^-_\Phi$: it consists
precisely of the objects whose bounded truncations vanish in the
algebraic quotient. In particular this quotient can have nonzero
objects with zero cohomology in every degree. Thus knowing its
heart is not enough to recover the quotient.

One might now hope that the difference disappears as soon as the
localized ring becomes artinian. It does not. For a minimal
nonmaximal prime $\mathfrak p$, the ring $R_{\mathfrak p}$ is
artinian, but the algebraic quotient at $R\setminus\mathfrak p$
still has a nonnoetherian Balmer spectrum. The criterion in
Corollary~\ref{thm:localized-noetherian-spectrum} instead requires
every prime surviving localization to be maximal in the original
ring. This is another distinction from the geometry of perfect
complexes: the localized coefficient ring does not by itself
determine the half-bounded algebraic quotient.

What if we remove a larger thick subcategory first? Could that
destroy the growth distinctions even though they survive scalar
localization? The question is meaningful for subcategories generated
by unbounded objects, where ordinary support gives little control
over uniform annihilation. We use the commutation of algebraic
localization with thick quotients to identify the objects being
killed, and check their images over the DVR. If those images are
uniformly annihilated, they belong to every ideal under consideration.
The inclusions then remain strict in the quotient.

\phantomsection\label{intro:relative-growth}
More precisely, under the hypotheses of
Theorem~\ref{intro:valuation-retraction} with $F$, let
$\mathcal C=\operatorname{Thick}(\mathcal E)$ be a tensor ideal and
$\Phi\subseteq\operatorname{End}(\mathbf1)$ multiplicative. Suppose
\[
 F(\phi)\ne0\quad(\phi\in\Phi),\qquad
 \forall E\in\mathcal E\ \exists a(E)\ge1:
             \pi^{a(E)}\operatorname{id}_{F(E)}=0.
\]
For $q:\mathcal K\to(\mathcal K/\mathcal C)_\Phi$, the ideals
$I_f=\sqrt{\operatorname{Thick}^{\otimes}(qK_f(s))}$ satisfy
\[
 I_f\subseteq I_g\ \Longleftrightarrow\ f\le_\mu g,
 \qquad I_f\vee I_g=I_{\max(f,g)},\qquad
 I_f\cap I_g=I_{\min(f,g)}.
\]
The quotient spectrum is therefore nonnoetherian
(Corollary~\ref{thm:central-relative-growth}). The exponent is allowed
to depend on the generator. The hypothesis concerns the identity
morphism, rather than only individual cohomology modules. This is
the condition to verify in the relative applications, including
Theorem~\ref{thm:central-relative-scheme} and
Example~\ref{ex:unbounded-relative-central}.

Are powers of degree-zero scalars the only way to obtain these
phenomena? Modular representations immediately lead us to a
different question. Over a field $k$ of characteristic $p\mid|G|$,
the degree-zero endomorphisms of the trivial representation form
$k$. There is no nonzero nonunit scalar to send to a positive power
of a DVR uniformizer. The preceding argument cannot start in that
form. But group cohomology supplies positive-degree classes, and a
cyclic subgroup of order $p$ has a degree-two periodicity class.
Can its powers take the place of the scalar powers?

They can, but the estimates must change. A degree-zero scalar
preserves degree, whereas a positive-degree class moves the Hom
group on which it acts. We estimate simultaneously how large a
power is needed to annihilate
$\operatorname{Hom}(X,\mathbf1[m])$ and how far it moves $m$.
The $t$-structure bounds the degrees that contribute. This gives
our second construction.

\begin{introthm}\label{intro:positive-growth}
Let $\mathcal T$ satisfy the coproduct hypotheses of
Theorem~\ref{intro:valuation-retraction}, and carry a $t$-structure
with $\mathcal T=\bigcup_b\mathcal T^{\le b}$,
$\mathbf1\in\mathcal T^\heartsuit$ and
$\mathcal T^{\le a}\otimes\mathcal T^{\le b}
\subseteq\mathcal T^{\le a+b}$.
Let $u:\mathbf1\to\mathbf1[q]$ be nonnilpotent, with $q>0$ even,
and assume $u^a_{\operatorname{Cone}(u^a)}=0$ for every $a\ge1$.
For $0<d<1$, set
\[
 E_d=\bigoplus_{i\ge0}
       \operatorname{Cone}(u^{\lfloor(i+1)^d\rfloor})[i-1].
\]
Then
\[
 \sqrt{\operatorname{Thick}^{\otimes}(E_e)}
       \subsetneq\sqrt{\operatorname{Thick}^{\otimes}(E_d)}
                      \qquad(0<e<d<1).
\]
Thus $\operatorname{Spc}(\mathcal T)$ is nonnoetherian.
\end{introthm}

Theorem~\ref{thm:abstract-periodic-growth} proves these estimates.
For a cyclic subgroup, the norm resolution supplies the class.
There is still an obstacle in passing to $G$: induction is not a
tensor functor, so it need not carry a strict tensor-generation
relation to a strict one. Proposition~\ref{prop:local-modular-growth}
uses restriction, the projection formula and Mackey decomposition
to control the additional summands. It proves strictness for all
tensor powers of the induced objects and shows that the chain
remains strict after quotienting by the tensor ideal of the regular
module. The norm resolution and Mackey formula are classical;
their use here must retain the growth estimates through both steps.

The two arguments give separate restrictions on noetherianity. The
coefficient ring can already contribute infinitely many radical
thick tensor ideals. A bad residue characteristic for $G$ contributes another
chain through group cohomology, even when the coefficient ring is
a field. In particular,
$\mathbf D^-_{\mathrm{fg}}(kG)$ can have the finite-length heart
$kG\text{-}\mathrm{mod}$ and a nonnoetherian Balmer spectrum.
Finiteness in each cohomological degree does not prevent the
unbounded growth of the constructions across degrees.

For a commutative noetherian ring $R$, we obtain
\[
 \boxed{\operatorname{Spc}(\mathbf D^-_{\mathrm{fg}}(R))
       \text{ is noetherian}\quad\Longleftrightarrow\quad
       R\text{ is artinian}.}
\]
The artinian implication is the Matsui--Takahashi result recalled
above. For the converse, a dominating DVR and the complexes constructed above
give the required spectral embedding
(Corollary~\ref{thm:noetherian-balmer-spectrum}). This identifies
positive Krull dimension as the ring-theoretic obstruction. Artinian
rings also provide a useful comparison: their half-bounded categories
still contain nondualizable objects, but their spectra are finite.
The theorem thus singles out a property of the coefficients, not
merely the presence of an infinite tail.

The affine criterion suggests that dimension should remain decisive
on schemes. Open restriction confirms this for arbitrary nonempty
noetherian schemes. For finite groups, restriction to the coefficient
ring and the cyclic-subgroup construction combine to give the group
criterion. For arithmetic spectra, extension to
$H\mathbb Z_{(p)}$ relates each noninvertible prime to the DVR
calculation. The results are summarized below; all spectra in the
second column are Balmer spectra.
\begin{center}
\renewcommand{\arraystretch}{1.18}
\begin{tabular}{@{}p{.19\textwidth}p{.77\textwidth}@{}}
Category & Noetherianity criterion\\\hline
$\mathbf D^-_{\mathrm{Coh}}(X)$ &
$\dim X=0$ for nonempty noetherian $X$
(Corollary~\ref{cor:scheme-noetherianity-criterion}).\\
$\mathbf D^-_{\mathrm{fg}}(RG)$ &
$R$ is artinian and $|G|\in R^\times$;
after the regular-module quotient and algebraic localization at
$\Phi$, precisely $|G|\in(\Phi^{-1}R)^\times$
(Theorem~\ref{thm:all-group-noetherian-ring}).\\
$\mathcal K_{S^{-1}\mathbb Z}$ &
$S^{-1}\mathbb Z=\mathbb Q$, where $\mathcal K$ consists of
bounded-below degreewise finite spectra
(Corollary~\ref{cor:arithmetic-spectra-noetherian}).
\end{tabular}
\end{center}
The relative scheme result retains the stated condition on the
identity morphisms after pullback to the chosen traits. In stable
homotopy theory, nonnoetherianity also follows from classical
chromatic theory. The
valuation construction gives the more specific DVR spectral
retractions and strict growth after quotienting by uniform torsion.
These stronger comparisons use the compatibility of the lattice
operations in Theorem~\ref{thm:valuation-growth-retraction}.

What should be the corresponding problem for bounded-below
complexes?
Over a general noetherian ring, bounded-below complexes need not
remain bounded below under derived tensor product. We therefore
use the action of $(\mathbf D^-_{\mathrm{fg}}(R))^{\mathrm{op}}$ given by
\[
 A\odot Y=\mathbf R\operatorname{Hom}_R(A,Y),\qquad
 A\in\mathbf D^-_{\mathrm{fg}}(R),
\]
on $\mathbf D^+_{\mathrm{fg}}(R\text{-}\mathrm{Mod})$ and on
$\mathbf D^+(R\text{-}\mathrm{mod})$. One category consists of derived objects with finite cohomology
modules, whereas the other is defined using bounded-below complexes
whose terms are finite modules. We construct
the Hom action and the localization comparisons in each category.
This lets us ask the same questions about perfect generation,
uniform annihilation and reconstruction from bounded truncations
without identifying the two positive categories.

We return at the end to what cohomology can recover. Postnikov
closure forgets the growth distinctions in our examples:
\[
 \operatorname{PC}^-(\mathcal J)
 =\mathbf D^-_{\mathrm{fg}}(R)_{\operatorname{Supp}\mathcal J}
 =\mathcal J^{\mathrm{tame}}.
\]
The last expression is the tame closure of Matsui--Takahashi
\cite[Lemma~5.8]{Matsui_2017}. Their support results identify the
common bounded and perfect parts. The chains constructed above show how many different radical ideals
can have those same parts and the same Postnikov closure. Thus the bounded part and the Postnikov closure
together do not determine a radical tensor ideal. We use these
comparisons to describe the additional information in the radical
ideals beyond the known support classification. The complete fourteen-condition statements,
with both positive versions, are in
Appendix~\ref{app:extended-criteria}.

The examples also leave a natural question for further work. How
generally is noetherianity of a half-bounded spectrum related to
the recovery of its ideals from bounded truncations? The ring,
scheme, representation and arithmetic cases give different tests
of this question. Conjecture~\ref{conj:postnikov-noetherianity}
formulates it under the hypotheses stated there.

Section~\ref{sec:scalar-growth} constructs the complexes from powers
of endomorphisms and proves the comparison with a DVR. Section~\ref{sec:relative-growth} treats
thick quotients and algebraic localization.
Section~\ref{sec:periodic-groups} gives the positive-degree argument
and its finite-group application, and
Section~\ref{sec:geometric-arithmetic} treats schemes and arithmetic
spectra. Section~\ref{sec:cohomological-information} studies the
cohomological comparisons and the concluding conjecture.
Appendices~\ref{app:central-localization}--\ref{app:extended-criteria}
contain the general quotient formulas, positive lifting proofs and
expanded criteria.

\begin{ack}
The author would like to thank Hongxing Chen and Xiaohu Chen
for discussions on the commutation of localization with Serre
and Verdier quotients in abstract tensor abelian and tensor
triangulated categories
(Propositions~\ref{thm:localization-commutes-serre-quotients}
and~\ref{thm:localization-commutes-triangulated-quotients}).
These discussions also led the author to study further the commutation
of geometric localization with quotients by thick subcategories in
$\mathbf D^+_{\mathrm{fg}}(R)$ and $\mathbf D^-_{\mathrm{fg}}(R)$.
\end{ack}

\subsection*{AI Declaration}
The framework of this paper, including its abstract and relative
formulations, was developed independently by the author in handwritten
notes before the AI-assisted discussions concerning this manuscript. The commutative-algebra
results in the author's original manuscript were obtained by the
author. Conjecture~\ref{conj:postnikov-noetherianity} was proposed
by the author.

The author used OpenAI ChatGPT Plus (model label: Sol) to assist with
organizing and rewriting the author's material, checking proof
details, searching the literature and verifying references.
The tool also assisted with identifying and developing the following
applications and working on their proofs:
\begin{itemize}[leftmargin=*,nosep]
\item Algebraic geometry: Corollary~\ref{cor:scheme-noetherianity-criterion},
Theorem~\ref{thm:central-relative-scheme},
and Examples~\ref{ex:unbounded-relative-central}
and~\ref{ex:projective-wide-postnikov}.
\item Stable homotopy theory: Corollary~\ref{cor:arithmetic-spectra-noetherian}
and the discussion following it.
\item Finite-group representations:
Proposition~\ref{prop:local-modular-growth},
Theorem~\ref{thm:all-group-noetherian-ring}, and
\hyperref[obs:group-postnikov-collapse]{the modular Postnikov observation}.
\end{itemize}
The author takes responsibility for the mathematical content,
references and final text.
\section{Tensor ideals and maps to discrete valuation rings}
\label{sec:scalar-growth}\label{sec:spectral-criteria}\label{sec:preliminaries}

The categories whose Balmer spectra occur below are essentially
small. All subcategories are full and closed under isomorphism.
We use cohomological grading, so $H^i(X[n])=H^{i+n}X$.
Unless stated otherwise, $R$ is a commutative noetherian ring,
$R\text{-}\mathrm{mod}$ denotes finite modules, and
$R\text{-}\mathrm{Mod}$ denotes all modules. We write
$\mathbf D^-_{\mathrm{fg}}(R)$ for bounded-above derived objects
with finitely generated cohomology and $\mathbf D^{\mathrm{perf}}(R)$ for perfect complexes.
For either sign,
\[
 \operatorname{Supp}X=\bigcup_i\operatorname{Supp}_R H^iX,
 \qquad
 \operatorname{Ann}_R X=\ker(R\longrightarrow\operatorname{End}_{\mathbf D(R)}X).
\]
The latter annihilates the identity, rather than only the cohomology.
For a triangle $A\to B\to C\to A[1]$,
$\operatorname{Ann}A\operatorname{Ann}C\subseteq\operatorname{Ann}B$
\cite[Lemma~2.5(1)]{Matsui_2017}.
We use the standard one-sided truncations and the Koszul complex
$K^\bullet(x_1,\ldots,x_r)$, with one-element factors in degrees $-1,0$.

The Balmer support is
\[
 \operatorname{Spp}_{\mathcal K}(X)
   =\{\mathcal P\in\operatorname{Spc}(\mathcal K):X\notin\mathcal P\}.
\]
Its closed sets are arbitrary intersections of such supports.
Write $\operatorname{Thick}^{\otimes}(\mathcal E)$ for the thick tensor
ideal generated by $\mathcal E$ and
\[
 \sqrt{\mathcal I}=\{X:X^{\otimes n}\in\mathcal I\text{ for some }n\ge1\}.
\]
In the derived categories of modules, coproducts with a common
one-sided cohomological bound and only finitely many nonzero summands
in each cohomological degree retain finitely generated cohomology.
This follows from exactness of module coproducts; these are the
coproducts used below.

\subsection{Powers of endomorphisms}
Let $\mathcal U$ be a symmetric tensor triangulated category with
countable coproducts, preserved by tensoring with any object. Let
$\mathcal K\subseteq\mathcal U$ be an essentially small full thick
tensor subcategory containing the unit $\mathbf1$ and the objects
\begin{equation}\label{eq:weighted-unit-sums}
 P_{\mathbf d}=\coprod_{n\ge0}\mathbf1[n]^{\oplus d_n}
 \qquad(d_n\in\mathbb N\text{ finite for every }n).
\end{equation}
All coproducts in the following construction are computed in
$\mathcal U$. Thick ideals will only be closed under finite constructions. Fix $s\in\operatorname{End}(\mathbf1)$, and
use the unit constraints to regard $s$ as a natural scalar on every
object. We write
\[
 C_a(s)=\operatorname{Cone}(s^a:\mathbf1\to\mathbf1),\qquad
 K_f(s)=\coprod_{n\ge0}C_{f(n)}(s)[n],\qquad
 K_0(s)=0,\quad K_\infty(s)=\mathbf1.
\]
Here $a\ge1$ and $f:\mathbb N\to\mathbb Z_{>0}$ is nondecreasing.
These objects belong to $\mathcal K$: each displayed sum of cones
is the cone of a diagonal endomorphism of a suitable $P_{\mathbf d}$.
We omit $s$ when it is fixed. Following \cite[Section~3]{sanders2026balmer}, write
\begin{equation}\label{eq:mu-growth-convention}
 f\le_\mu g\quad\Longleftrightarrow\quad
 f(n)\le A g(2^kn)\ (n\gg0)
 \quad\text{for some }A\in\mathbb Z_{>0},\ k\in\mathbb N.
\end{equation}
Write $f\sim_\mu g$ for the associated equivalence relation, and
$L_\mu$ for the lattice of positive classes with a bottom $0$ and
a top $\infty$ adjoined. Joins and meets of positive classes are
represented by pointwise maximum and minimum. This is the lattice
$(\mathsf{ASeq}/\mu)_+$ of
\cite[Section~3 and Theorem~11.1]{sanders2026balmer}.

For $X\in\mathcal K$, abbreviate
\[
 \beta_{\mathcal K}(X)
   =\sqrt{\operatorname{Thick}^{\otimes}(X)}.
\]
The principal radical ideals form a bounded distributive lattice
$\mathsf{PRad}(\mathcal K)$, with bottom $\beta_{\mathcal K}(0)$,
top $\mathcal K$, and operations
\begin{equation}\label{eq:principal-radical-product}
 \beta(X)\vee\beta(Y)=\beta(X\oplus Y),\qquad
 \beta(X)\cap\beta(Y)=\beta(X\otimes Y).
\end{equation}
The meet formula follows by testing containment in primes
\cite[Lemma~4.2]{balmer2005spectrum}. The lattice description is
also given in \cite[Section~3]{kockpitsch2017hochster}.
In a non-rigid category $\beta(0)$ need not be the zero ideal.
Write $\mathsf{Rad}(\mathcal K)$ for all radical thick tensor ideals,
ordered by inclusion.

\begin{lemma}\label{lem:scalar-koszul-calculus}
In any tensor triangulated category,
\[
 s^{2a}\operatorname{id}_{C_a}=0,\qquad
 C_a\otimes C_b\simeq C_a\oplus C_a[1]\quad(b\ge2a).
\]
There is also a triangle
$C_a\to C_{a+b}\to C_b\to C_a[1]$ for $a,b\ge1$.
\end{lemma}
\begin{proof}
Set $u=s^a$ and write its triangle as
$\mathbf1\xrightarrow{u}\mathbf1\xrightarrow{i}C_a
\xrightarrow{p}\mathbf1[1]$.
Naturality gives $u_{C_a}i=iu=0$, so $u_{C_a}=hp$ for some $h$.
Also $pu_{C_a}=u[1]p=0$. Hence $u_{C_a}^2=0$.
If $b\ge2a$, tensoring the triangle defining $C_b$ with $C_a$
gives the split cone of $s^b\operatorname{id}_{C_a}=0$.
The last triangle is the octahedral triangle for
$\mathbf1\xrightarrow{s^a}\mathbf1\xrightarrow{s^b}\mathbf1$.
\end{proof}

The proofs of \cite[Lemmas~7.20--7.21 and~10.8;
Propositions~10.12 and~10.15]{sanders2026balmer} use diagonal
cones, tensor powers and finite multiplicities in each degree.
We apply these constructions to the objects $C_a(s)$ in $\mathcal K$.

\begin{lemma}
\label{lem:koszul-growth-bound}
In the above setting, let $g:\mathbb N\to\mathbb Z_{>0}$ be
nondecreasing and let
\[
 X=\coprod_{n\ge0}\coprod_{j=1}^{m_n}C_{a_{nj}}[n],
 \qquad m_n<\infty,\quad a_{nj}\ge1.
\]
Suppose $a_{nj}\le A g(2^kn)$ for all $j$ and $n\gg0$, with
$A\ge1$ and $k\ge0$ integers. After increasing $A$ to cover the
initial degrees, one has, for $q=2^{k+1}$,
\[
 X^{\otimes q}\in\operatorname{Thick}^{\otimes}(K_g).
\]
In particular $X\in\beta_{\mathcal K}(K_g)$.
\end{lemma}
\begin{proof}
The finitely many initial exponents can be absorbed into $A$
because $g\ge1$. Expand $X^{\otimes q}$ in $\mathcal U$. Index
its blocks by tuples $\tau=((n_1,j_1),\ldots,(n_q,j_q))$ and set
$N=\sum_\nu n_\nu$. Choose a factor of least degree
$n_{\nu_0}\le\lfloor N/q\rfloor$, and denote its exponent by
$b_\tau=a_{n_{\nu_0}j_{\nu_0}}$. Then
\[
 b_\tau\le A g(2^k\lfloor N/q\rfloor)
            \le A g(\lfloor N/2\rfloor).
\]
There are finitely many tuples of each total degree $N$; call their
number $d_N$. First construct only the selected factors:
\[
 Y_0=\coprod_\tau C_{b_\tau}[N].
\]
Let $e_0=d_0$, $e_j=\max(d_{2j-1},d_{2j})$ for $j\ge1$, and
$Q=P_{\mathbf e}$. In
\[
 K_g\otimes Q
   \simeq\coprod_{i,j\ge0}C_{g(i)}[i+j]^{\oplus e_j}
\]
choose $d_N$ copies at $(i,j)=(\lfloor N/2\rfloor,\lceil N/2\rceil)$.
They give the direct summand
\[
 Z_1=\coprod_{N\ge0}
             C_{g(\lfloor N/2\rfloor)}[N]^{\oplus d_N}
 \in\operatorname{Thick}^{\otimes}(K_g).
\]
Let $Z_t$ have the same blocks with exponents multiplied by $t$.
The coproducts of the octahedral triangles in
Lemma~\ref{lem:scalar-koszul-calculus} give triangles
$Z_{t-1}\to Z_t\to Z_1\to Z_{t-1}[1]$.
Thus $Z_{2A}$ is obtained by $2A-1$ extensions.
On the block indexed by $\tau$, let $d:Z_{2A}\to Z_{2A}$ be
multiplication by $s^{b_\tau}$. Since
$a=2A g(\lfloor N/2\rfloor)\ge2b_\tau$, symmetry and
Lemma~\ref{lem:scalar-koszul-calculus} give on each block
\[
 \operatorname{Cone}(s^{b_\tau}:C_a\to C_a)
 \simeq C_a\otimes C_{b_\tau}
 \simeq\operatorname{Cone}(s^a:C_{b_\tau}\to C_{b_\tau})
 \simeq C_{b_\tau}\oplus C_{b_\tau}[1].
\]
Thus $\operatorname{Cone}(d)\simeq Y_0\oplus Y_0[1]$.
Hence $Y_0\in\operatorname{Thick}^{\otimes}(K_g)$.

For each tuple, order the other $q-1$ exponents arbitrarily.
Taking on $Y_0$ the diagonal scalar endomorphism given by the first
of these exponents, and then taking its cone, inserts the next
factor $C_{a_{n_\nu j_\nu}}$ in every block simultaneously.
Repeat for the other exponents. After $q-1$ cones the result is
exactly $X^{\otimes q}$, up to the symmetry constraints.
Every intermediate coproduct is in $\mathcal K$, since it is a
cone of an endomorphism of the preceding object. This construction
uses one tensor factor $Q$, a direct summand, and finitely many
extensions and cones; it does not close an ideal under coproducts.
\end{proof}

\subsection{Comparison with a discrete valuation ring}

\begin{theorem}
\label{thm:valuation-growth-retraction}
Let $\mathcal K\subseteq\mathcal U$ and $s$ satisfy the assumptions
above.
\begin{enumerate}
\item The assignments
\[
 \gamma_s([f]_\mu)=\beta_{\mathcal K}(K_f(s)),\qquad
 \gamma_s(0)=\beta_{\mathcal K}(0),\qquad
 \gamma_s(\infty)=\mathcal K
\]
define a homomorphism of bounded distributive lattices
$L_\mu\to\mathsf{PRad}(\mathcal K)$.
\item Let $\mathcal L$ satisfy the same assumptions in an ambient
category, and let $F:\mathcal K\to\mathcal L$ be an exact symmetric
monoidal functor preserving the coproduct
$P_{\mathbf1}=\coprod_{n\ge0}\mathbf1[n]$, including its coprojections.
Set $t=F(s)$. If
$\gamma_t:L_\mu\xrightarrow{\sim}\mathsf{PRad}(\mathcal L)$
is an isomorphism, then $\gamma_s$ is injective, and
\[
 \beta_{\mathcal K}(K_f(s))\subseteq\beta_{\mathcal K}(K_g(s))
       \quad\Longleftrightarrow\quad f\le_\mu g.
\]
For $F_*(\mathcal I)=
\sqrt{\operatorname{Thick}^{\otimes}(F(\mathcal I))}$,
$\sigma_s=\gamma_s\gamma_t^{-1}$ is a lattice section of $F_*$
on principal radical ideals.
\item Consequently, this section extends to
\[
 \Sigma_s:\mathsf{Rad}(\mathcal L)\longrightarrow\mathsf{Rad}(\mathcal K),
 \qquad
 \Sigma_s(\mathcal I)=
   \bigvee_{\beta_{\mathcal L}(E)\subseteq\mathcal I}
                    \sigma_s(\beta_{\mathcal L}(E)),
\]
preserving arbitrary joins and finite meets, with
$F_*\Sigma_s=\operatorname{id}$.
The spectral map $\theta=\operatorname{Spc}(F)$ has a spectral
retraction
\begin{equation}\label{eq:valuation-spectral-retraction}
 \varrho_s:\operatorname{Spc}(\mathcal K)
          \longrightarrow\operatorname{Spc}(\mathcal L),
 \qquad \varrho_s\theta=\operatorname{id}.
\end{equation}
In particular, $\theta$ is a topological embedding.
\end{enumerate}
\end{theorem}
\begin{proof}
Lemma~\ref{lem:koszul-growth-bound} shows that $f\le_\mu g$ implies
$\beta(K_f)\subseteq\beta(K_g)$. Thus $\gamma_s$ is well defined
and order preserving. Summing octahedral triangles in $\mathcal U$
gives $K_f\to K_{f+g}\to K_g\to K_f[1]$. Since
$f,g\le_\mu f+g\sim_\mu\max(f,g)$, this proves the join formula.

For the meet, write $h(n)=\min(f(n),g(n))$ and expand
$K_f\otimes K_g$ into blocks
$(C_{f(i)}\otimes C_{g(j)})[i+j]$. Let
$b_{ij}=\min(f(i),g(j))$ and $c_{ij}=\max(f(i),g(j))$.
As $b_{ij}\le h(i+j)$, Lemma~\ref{lem:koszul-growth-bound} gives
\[
 Y=\coprod_{i,j\ge0}C_{b_{ij}}[i+j]\in\beta(K_h).
\]
The cone of the diagonal endomorphism of $Y$ with $(i,j)$ block
$s^{c_{ij}}$ is $K_f\otimes K_g$. Hence
$\beta(K_f\otimes K_g)\subseteq\beta(K_h)$.
The reverse containment follows from $h\le_\mu f,g$ and
\eqref{eq:principal-radical-product}. This proves (1), including
the bottom and top values.

Preservation of $P_{\mathbf1}$ identifies the image of its diagonal
endomorphism $\operatorname{diag}(s^{f(n)})$ with
$\operatorname{diag}(t^{f(n)})$. Exactness therefore gives
$F(K_f(s))\simeq K_f(t)$, and consequently
\[
 F_*\gamma_s=\gamma_t.
\]
If $\gamma_t$ is an isomorphism, this proves the injectivity, order
reflection and section asserted in (2).

For (3), radical ideals form the ideal completion of the lattice
of principal radical ideals. Thus $\sigma_s$ extends by the displayed
join formula to a frame section $\Sigma_s$ of $F_*$. Hochster duality
gives its spectral retraction
\cite[Section~1.2, Theorem~3.1.9 and Corollary~3.4.2]{kockpitsch2017hochster}:
\[
 \varrho_s(\mathcal P)=
 \{E\in\mathcal L:\sigma_s(\beta_{\mathcal L}(E))
                                      \subseteq\mathcal P\}.
\]
For $\beta_{\mathcal L}(E)=\gamma_t(f)$, including $f=0,\infty$,
\[
 \varrho_s^{-1}(\operatorname{Spp}_{\mathcal L}(E))
       =\operatorname{Spp}_{\mathcal K}(K_f(s)),\qquad
 \varrho_s\operatorname{Spc}(F)=\operatorname{id}.
\]
\end{proof}

\smallskip\noindent\emph{Application to rings.}
Take $\mathcal U=\mathbf D(R)$ and
$\mathcal K=\mathbf D^-_{\mathrm{fg}}(R)$.
The objects $P_{\mathbf d}$ are bounded-above complexes of finite
free modules. For a map $u:R\to V$ to a DVR such that
$u(r)=\varepsilon\pi^c$, with $\varepsilon\in V^\times$ and $c\ge1$,
extension of scalars $F=-\otimes_R^{\mathbf L}V$ preserves the
required coproduct and satisfies
\[
 F(K_f(r))\simeq K_f(u(r))
       \simeq T_{cf},\qquad
 T_h=\bigoplus_{n\ge0}V/\pi^{h(n)}V[n].
\]
By \cite[Theorems~10.14 and~11.1]{sanders2026balmer} and
$cf\sim_\mu f$, the target map $\gamma_{u(r)}$ is an isomorphism.
Consequently, with $J_f=\beta_{\mathbf D^-_{\mathrm{fg}}(R)}(K_f(r))$,
\begin{equation}\label{eq:exact-koszul-growth-comparison}
 J_f\subseteq J_g\quad\Longleftrightarrow\quad f\le_\mu g,
 \qquad
 J_f\vee J_g=J_{\max(f,g)},\quad J_f\cap J_g=J_{\min(f,g)}.
\end{equation}
Here $\beta(0)=0$. We write $\Sigma_r$ and $\varrho_r$ for the
resulting section and spectral retraction.

\begin{corollary}\label{prop:dvr-radical-lifting}
Let $u:R\to V$ be a homomorphism to a DVR with
$u(r)=\varepsilon\pi^c$, $\varepsilon\in V^\times$ and $c>0$.
For $F=-\otimes_R^{\mathbf L}V$, every
$E\in\mathbf D^-_{\mathrm{fg}}(V)$ admits
$X\in\mathbf D^-_{\mathrm{fg}}(R)$ such that
\begin{equation}\label{eq:principal-radical-lifting}
 \sqrt{\operatorname{Thick}^{\otimes}(E)}
       =\sqrt{\operatorname{Thick}^{\otimes}(FX)}.
\end{equation}
The map $\operatorname{Spc}(F)$ is a topological embedding and has
a spectral retraction.
\end{corollary}
\begin{proof}
By \cite[Theorems~10.14 and~11.1]{sanders2026balmer}, every principal
radical ideal over $V$ is represented by $0$, $V$, or
$T_f=\bigoplus_{n\ge0}V/\pi^{f(n)}V[n]$ with $f$ positive and
nondecreasing. Take $X=0$, $R$, or $K_f(r)$ respectively. Then
$F(K_f(r))\simeq T_{cf}$ and $cf\sim_\mu f$, giving the equality.
For $\theta=\operatorname{Spc}(F)$, it follows that
\[
 \operatorname{Spp}_V(E)
       =\operatorname{Spp}_V(FX)
       =\theta^{-1}(\operatorname{Spp}_R(X)).
\]
These supports distinguish points and generate all closed subsets,
so $\theta$ is a topological embedding. The spectral retraction
follows from Theorem~\ref{thm:valuation-growth-retraction}.
\end{proof}

\begin{corollary}
\label{thm:noetherian-balmer-spectrum}
Let $R$ be a commutative noetherian ring. Then
\[
 R\text{ is artinian}
 \quad\Longleftrightarrow\quad
 \operatorname{Spc}(\mathbf D^-_{\mathrm{fg}}(R))
           \text{ is a noetherian topological space}.
\]
\end{corollary}
\begin{proof}
If $R$ is artinian, \cite[Theorem~6.5]{Matsui_2017} identifies the
Balmer spectrum with the finite space $\operatorname{Spec}R$.
Otherwise choose $\mathfrak p\subsetneq\mathfrak m$, with
$\mathfrak m$ maximal, and a DVR $V$ dominating the nonfield
local domain $(R/\mathfrak p)_{\mathfrak m/\mathfrak p}$
\cite[Tags~00KH and~00PH]{stacks-project}. For
$r\in\mathfrak m\setminus\mathfrak p$, the composite $R\to V$
sends $r$ to a nonzero nonunit. Corollary~\ref{prop:dvr-radical-lifting} gives an embedding
$\operatorname{Spc}(\mathbf D^-_{\mathrm{fg}}(V))\hookrightarrow
\operatorname{Spc}(\mathbf D^-_{\mathrm{fg}}(R))$.
The former is nonnoetherian
\cite[Remark~12.12]{sanders2026balmer}, whereas every subspace of
a noetherian space is noetherian. This proves the converse.
\end{proof}
\section{Growth after quotients and algebraic localization}
\label{sec:relative-growth}

The kernel of central localization consists of objects annihilated
by some scalar in the multiplicative set. Localization of each
cohomology module can have a larger kernel. We first record the
quotient formula.

Let $\mathcal T$ be an essentially small tensor triangulated category,
$C=\operatorname{End}_{\mathcal T}(\mathbf1)$ and
$\Phi\subseteq C$ multiplicative. Write
\[
 \mathcal S_\Phi=\{X\in\mathcal T:
       \phi\operatorname{id}_X=0\text{ for some }\phi\in\Phi\},
 \qquad\mathcal T_\Phi=\mathcal T/\mathcal S_\Phi.
\]

\begin{proposition}\label{prop:minimal-central-localization}
For $\pi:\mathcal T\to\mathcal T_\Phi$,
\begin{align}
 \mathcal S_\Phi
   &=\operatorname{Thick}^{\otimes}
       \{\operatorname{Cone}(\phi:\mathbf1\to\mathbf1):\phi\in\Phi\},
       \label{eq:early-central-cone-generators}\\
 \operatorname{Hom}_{\mathcal T_\Phi}(\pi X,\pi Y)
   &\cong\Phi^{-1}\operatorname{Hom}_{\mathcal T}(X,Y).
       \label{eq:early-central-hom}
\end{align}
For every thick subcategory $\mathcal C\subseteq\mathcal T$,
\begin{equation}\label{eq:minimal-central-relative-quotient}
 (\mathcal T/\mathcal C)_\Phi
   \simeq\mathcal T/\operatorname{Thick}(\mathcal C\cup\mathcal S_\Phi).
\end{equation}
If $\mathcal C$ is a tensor ideal, these are tensor functors.
\end{proposition}
\begin{proof}
For a triangle $X\to Y\to Z\to X[1]$, scalar naturality gives
$\operatorname{Ann}_C(X)\operatorname{Ann}_C(Z)
\subseteq\operatorname{Ann}_C(Y)$; together with shifts, retracts
and tensor products, this makes $\mathcal S_\Phi$ a thick tensor ideal.
The cone of $\phi:\mathbf1\to\mathbf1$ is killed by $\phi^2$,
while $\phi_X=0$ makes $X$
a retract of $\operatorname{Cone}(\phi)\otimes X$, proving the first
formula. The Hom formula is central localization
\cite[Corollary~3.10]{balmer2010spectra}.
If $q(\phi_X)=0$ in $\mathcal T/\mathcal C$, then
\[
 q\operatorname{Cone}(\phi_X)\simeq qX\oplus qX[1].
\]
Thus the uniform torsion objects in that quotient are
$\operatorname{smd}(q\mathcal S_\Phi)$, and the Verdier universal
property gives \eqref{eq:minimal-central-relative-quotient}
\cite[Tag~05RJ]{stacks-project}. The general module statement and
image formulas are in Proposition~\ref{prop:triangulated-localization-hom}.
\end{proof}

\begin{example}\label{ex:triangulated-unbounded-torsion}
For $R=k[t]$,
\[
 U=\bigoplus_{n\ge1}R/(t^n)[-n]\in\mathbf D^+_{\mathrm{fg}}(R)
\]
satisfies $U_t=0$, but $t^m\operatorname{id}_U\ne0$ for every
$m\ge0$, as seen on $H^{m+1}U=R/(t^{m+1})$. Thus geometric
localization kills $U$, whereas the central quotient retains it.
For bounded complexes the kernels agree by the annihilator
product rule and a finite Postnikov filtration.
\end{example}

\subsection{Passing to a thick quotient}

\begin{corollary}
\label{thm:central-relative-growth}
Retain the hypotheses of
Theorem~\ref{thm:valuation-growth-retraction}, with
$F:\mathcal T\to\mathbf D^-_{\mathrm{fg}}(V)$ and
$F(s)=\varepsilon\pi^c$, $c>0$.
Let $\mathcal C=\operatorname{Thick}(\mathcal E)$ be a tensor ideal
and let $\Phi\subseteq\operatorname{End}_{\mathcal T}(\mathbf1)$
be multiplicative. Suppose
\[
 F(\phi)\ne0\quad(\phi\in\Phi),\qquad
 \forall E\in\mathcal E\ \exists a\ge1:\quad\pi^a\operatorname{id}_{F(E)}=0.
\]
For the composite $q:\mathcal T\to(\mathcal T/\mathcal C)_\Phi$
and a positive nondecreasing $f$, let
\[
 \mathcal I_f=
 \sqrt{\operatorname{Thick}^{\otimes}(qK_f(s))}.
\]
Write $\mathcal S_\pi=\{Z:\pi^a\operatorname{id}_Z=0\text{ for some }a\ge1\}$. Then
\begin{gather}
 \mathcal I_f\subseteq\mathcal I_g
       \quad\Longleftrightarrow\quad f\le_\mu g,
       \label{eq:central-relative-exact-growth}\\
 \mathcal I_f\vee\mathcal I_g=\mathcal I_{\max(f,g)},\qquad
 \mathcal I_f\cap\mathcal I_g=\mathcal I_{\min(f,g)}.
       \label{eq:central-relative-growth-lattice}
\end{gather}
In particular $\operatorname{Spc}((\mathcal T/\mathcal C)_\Phi)$
is not noetherian.
\end{corollary}
\begin{proof}
Write $\beta$ for principal radicals and $F_*$ for radical ideal
extension.
By Proposition~\ref{prop:minimal-central-localization}, the target is
$\mathcal T/\mathcal D$, with
$\mathcal D=\operatorname{Thick}(\mathcal C\cup\mathcal S_\Phi)$.
Every nonzero $F(\phi)$ is a unit times a power of $\pi$. The
annihilator hypotheses and closure under finite constructions give
\[
 F(\mathcal D)\subseteq\mathcal S_\pi
       \subseteq\beta(K_{cf}(\pi))\qquad(f\ge1).
\]
Set $d=\sqrt{\mathcal D}$ and $b_f=\beta(K_f(s))$.
The radical-ideal correspondence for Verdier quotients
\cite[Proposition~3.11]{balmer2005spectrum} gives
\[
 q^{-1}\mathcal I_f=d\vee b_f,\qquad
 F_*(d\vee b_f)=\beta(K_{cf}(\pi)),
\]
\[
 \begin{aligned}
 f\le_\mu g&\Longrightarrow d\vee b_f\subseteq d\vee b_g\\
           &\Longrightarrow\beta(K_{cf}(\pi))\subseteq\beta(K_{cg}(\pi))\\
           &\Longrightarrow f\le_\mu g.
 \end{aligned}
\]
The first implication follows from
Theorem~\ref{thm:valuation-growth-retraction}(1), and the last from the
DVR classification. The same theorem gives the join formula, and
distributivity \cite[Theorem~3.1.9]{kockpitsch2017hochster} gives
\[
 (d\vee b_f)\cap(d\vee b_g)
       =d\vee(b_f\cap b_g)=d\vee b_{\min(f,g)}.
\]
For $f_j(n)=\lceil(n+1)^{1/j}\rceil$, $j\ge1$, the supports of
$qK_{f_j}(s)$ therefore form a strictly descending chain of closed
subsets.
\end{proof}
The annihilator exponent may depend on the generator $E$.
Here $f$ takes positive values. The object $K_1(s)$ may become zero
in the quotient, so its ideal can be the least ideal in this family.

\subsection{The spectrum of a central quotient}
For either half-bounded finite-cohomology category over $R$, write
$\mathcal S^\pm_\Phi$ for uniform scalar torsion as above. A subscript
$W$ denotes objects with cohomological support in $W$.

\begin{corollary}\label{thm:localized-noetherian-spectrum}
Let $\Phi\subseteq R$ be multiplicative, with
$U_\Phi=\{\mathfrak p\in\operatorname{Spec}R:\mathfrak p\cap\Phi=\varnothing\}$.
The following conditions are equivalent:
\begin{enumerate}
\item $\operatorname{Spc}(\mathbf D^-_{\mathrm{fg}}(R)/\mathcal S^-_\Phi)$
      is noetherian;
\item $U_\Phi\subseteq\operatorname{Max}R$;
\item $\Phi^{-1}R$ is artinian and
\[
 \overline L^-_\Phi:\mathbf D^-_{\mathrm{fg}}(R)/\mathcal S^-_\Phi
             \longrightarrow\mathbf D^-_{\mathrm{fg}}(\Phi^{-1}R)
\]
is an equivalence.
\end{enumerate}
In this case $U_\Phi$ is finite and clopen, and the localization map
identifies $\Phi^{-1}R$ with an artinian direct factor $eR$ of $R$.
For both signs,
\[
 \mathcal S^\pm_\Phi=\mathcal N^\pm_\Phi
 =\mathbf D^\pm_{\mathrm{fg}}(R)_{\operatorname{Spec}R\setminus U_\Phi},
 \qquad \overline L^\pm_\Phi\text{ is an equivalence},
\]
where $\mathcal N^\pm_\Phi=\{X:\Phi^{-1}H^iX=0\ (\forall i)\}$.
\end{corollary}
\begin{proof}
If $\mathfrak p\in U_\Phi$ is not maximal, choose a maximal
$\mathfrak m\supsetneq\mathfrak p$ and a DVR $V$ dominating
$(R/\mathfrak p)_{\mathfrak m/\mathfrak p}$, as in
Corollary~\ref{thm:noetherian-balmer-spectrum}. For the composite
$u:R\to V$, every $u(\phi)$, $\phi\in\Phi$, is nonzero, while some
$r\in\mathfrak m\setminus\mathfrak p$ has nonzero nonunit image.
Corollary~\ref{thm:central-relative-growth}, with $\mathcal C=0$,
therefore makes the central quotient spectrum nonnoetherian.
This proves $(1)\Rightarrow(2)$.

Assume (2). Generalization closure gives
$U_\Phi\subseteq\operatorname{Min}R\cap\operatorname{Max}R$.
Such primes are clopen components, and there are finitely many.
For the corresponding idempotent $e$,
\[
 \Phi^{-1}R\cong eR,\qquad \dim eR=0,\qquad
 \phi_0(1-e)=0\text{ for some }\phi_0\in\Phi.
\]
Every $\phi\in\Phi$ acts invertibly on the $eR$-factor, and
$\phi_0$ kills the complementary factor. The product decomposition
of either derived category proves $\mathcal S^\pm_\Phi=\mathcal N^\pm_\Phi$
and all the additional claims. Thus $(2)\Rightarrow(3)$.
Finally $(3)\Rightarrow(1)$ is
Corollary~\ref{thm:noetherian-balmer-spectrum} for $\Phi^{-1}R$.
\end{proof}

For $\Phi=R\setminus\mathfrak p$, this gives
\[
 \operatorname{Spc}(\mathbf D^-_{\mathrm{fg}}(R)/\mathcal S^-_\Phi)
   \text{ noetherian}
 \quad\Longleftrightarrow\quad
 \mathfrak p\in\operatorname{Min}R\cap\operatorname{Max}R.
\]
By contrast, the geometric target
$\mathbf D^-_{\mathrm{fg}}(R_{\mathfrak p})$ has noetherian spectrum
exactly when $\mathfrak p$ is minimal. A minimal nonmaximal prime
therefore separates the two noetherianity conditions.
\section{Periodic growth and modular representations}\label{sec:periodic-groups}

The second construction uses a positive-degree class on the tensor unit.

\begin{theorem}
\label{thm:abstract-periodic-growth}
Let $\mathcal T\subseteq\mathcal U$ satisfy the assumptions preceding
\eqref{eq:weighted-unit-sums}; in particular, $\mathcal T$ contains
each coproduct $P_{\mathbf d}$ in that formula.
Suppose $\mathcal T$ has a $t$-structure with
$\mathcal T=\mathcal T^-$, $\mathbf1\in\mathcal T^\heartsuit$, and
\[
 \mathcal T^{\le a}\otimes\mathcal T^{\le b}
                  \subseteq\mathcal T^{\le a+b}.
\]
Let $u:\mathbf1\to\mathbf1[q]$ be nonnilpotent, with $q>0$ even.
Let $P_a=\operatorname{Cone}(u^a)[-1]$ and assume
$u^a_{P_a}=0$ for every $a\ge1$.
For $0<d<1$ set
\[
 a_i(d)=\lfloor(i+1)^d\rfloor,\qquad
 E_d=\bigoplus_{i\ge0}P_{a_i(d)}[i].
\]
Then these objects belong to $\mathcal T$, and
\[
 \sqrt{\operatorname{Thick}^{\otimes}(E_e)}
   \subsetneq\sqrt{\operatorname{Thick}^{\otimes}(E_d)}
                       \qquad(0<e<d<1).
\]
Consequently $\operatorname{Spc}(\mathcal T)$ is not noetherian.
\end{theorem}
\begin{proof}
All sums used below exist in $\mathcal U$ and lie in $\mathcal T$:
each is a cone between coproducts of shifts of the unit with finite
multiplicity in each degree, or a retract of
a displayed tensor product. Existing coproducts of objects of
$\mathcal T^{\le B}$ remain in $\mathcal T^{\le B}$ whenever their
coproduct lies in $\mathcal T$, by orthogonality.
In particular $P_a\in\mathcal T^{\le0}$.
For $a\ge b$, tensoring the defining triangle and using $u^a_{P_b}=0$
gives
\begin{equation}\label{eq:periodic-cone-product}
 P_a\otimes P_b\simeq P_b\oplus P_b[qa-1].
\end{equation}

Write $\ell_X(n)$ for the least $N\ge0$ such that
$u^N\operatorname{Hom}(X,\mathbf1[m])=0$ for every $m\le n$,
allowing $\infty$. For $X\le B$ these Hom groups vanish when
$m<-B$. Naturality of $u$ makes the cohomological long exact
sequence $u$-linear. For a triangle $X\to Y\to Z\to X[1]$ it gives
\[
 \ell_Z(n)\le\ell_Y(n)+\ell_X(n+q\ell_Y(n)+1)
\]
when the right side is finite: after applying $u^{\ell_Y(n)}$,
each class comes from degree $m+q\ell_Y(n)-1$ for $X$.
Thus the objects with $\ell_X(n)=O((n+1)^e)$ form a thick class
$\mathcal C_e$ for $0<e<1$, also closed under shifts and retracts.

If $Y\le B$, the $i$-th summand of $E_e\otimes Y$ lies in
$\mathcal T^{\le B-i}$. Its Hom to $\mathbf1[m]$ vanishes for
$i>m+B$ by $t$-orthogonality. Only finitely many summands contribute
in each degree, and the whole $i$-th summand is annihilated by
$u^{a_i(e)}$. Hence
\[
 \ell_{E_e\otimes Y}(n)
       \le\max_{0\le i\le n+B}a_i(e)=O((n+1)^e),\qquad
 \operatorname{Thick}^{\otimes}(E_e)\subseteq\mathcal C_e.
\]
The second assertion uses the ordinary thick closure of the objects
$E_e\otimes Y$, not closure under arbitrary coproducts.

The image of $1$ in $\operatorname{Hom}(P_a,\mathbf1)$ has
$u$-annihilator exponent exactly $a$. Indeed for $0\le N<a$,
the triangle for $P_a$ makes
\[
 \operatorname{Hom}(\mathbf1,\mathbf1[qN])
       \longrightarrow\operatorname{Hom}(P_a,\mathbf1[qN])
\]
injective, since $\operatorname{Hom}(\mathbf1,\mathbf1[q(N-a)])=0$;
the image of $u^N$ is therefore nonzero. The equality $u^a_{P_a}=0$ gives the upper bound $a$. By \eqref{eq:periodic-cone-product},
$P_{a_i(d)}[ri]$ is a retract of $E_d^{\otimes r}$ for every $r\ge1$.
Thus $\ell_{E_d^{\otimes r}}(ri)\ge a_i(d)$, excluding all these
powers from $\mathcal C_e$ when $d>e$.

For the reverse inclusion, write $a_i=a_i(d)$ and $b_i=a_i(e)\le a_i$.
The object $T=\bigoplus_iP_{a_i}[i+qb_i]$ is a retract of
$E_d\otimes Q$, where $Q=\bigoplus_{j\ge0}\mathbf1[j]$:
retain the blocks $(i,j=qb_i)$. The diagonal morphism
$E_d\to T$ with components $u^{b_i}$ has cone
\[
 E_e[1]\oplus\bigoplus_iP_{b_i}[i+qa_i]
\]
by \eqref{eq:periodic-cone-product}. Its summand $E_e[1]$ therefore
belongs to $\operatorname{Thick}^{\otimes}(E_d)$.
Since no tensor power of $E_d$ belongs to $\mathcal C_e$ when $d>e$,
the radical inclusions are strict.
The supports for $d=1/2,1/3,1/4,\ldots$ consequently give a strictly
descending chain of closed subsets, by the radical classification
\cite[Theorem~4.10]{balmer2005spectrum}.
\end{proof}

\subsection{Transfer from a cyclic subgroup}

\begin{proposition}\label{prop:local-modular-growth}
Let $(R,\mathfrak m,k)$ be a nonzero commutative noetherian local
ring with $\operatorname{char}k=p>0$, and let $G$ be finite with
$p\mid |G|$. Choose $H=C_p\le G$, and let
$u\in\operatorname{Ext}^2_{RH}(R,R)$ be the class represented by the
periodic norm resolution. Let
\[
 P_a=\operatorname{Cone}(u^a)[-1],\qquad
 X_d=\operatorname{Ind}_H^G
       \left(\bigoplus_{i\ge0}P_{\lfloor(i+1)^d\rfloor}[i]\right)
       \quad(0<d<1).
\]
Then $X_d\in\mathbf D^-_{\mathrm{fg}}(RG)$,
\[
 RG\in\operatorname{Thick}^{\otimes}(X_d),\qquad
 H^0X_d\cong\operatorname{Ind}_H^G R,
\]
and
\[
 \sqrt{\operatorname{Thick}^{\otimes}(X_e)}
  \subsetneq\sqrt{\operatorname{Thick}^{\otimes}(X_d)}
  \qquad(0<e<d<1).
\]
Consequently
$\operatorname{Spc}(\mathbf D^-_{\mathrm{fg}}(RG)/
\operatorname{Thick}^{\otimes}(RG))$ is not noetherian.
\end{proposition}
\begin{proof}
Choose $H=C_p\le G$ and set $A=RH$.
The usual norm resolution supplies a periodic nonnilpotent class
$u\in H^2(H,R)$, with $H^{2a}(H,R)=R/pR$ for $a\ge1$
\cite[Examples~1.3 and~3.2]{brown2010lectures}.
Its $a$-th power is represented by the $2a$-extension whose middle
terms are finite free $A$-modules. Thus
\[
 P_a=\operatorname{Cone}(u^a:R\to R[2a])[-1]
\]
has a finite free representative in degrees $-2a+1,\ldots,0$.
Consequently $\operatorname{Hom}(P_a,P_a[2a])=0$, so $u^a_{P_a}=0$.
The tensor product is right $t$-exact, and the coproducts
$P_{\mathbf d}$ have finitely many summands in each degree.
Theorem~\ref{thm:abstract-periodic-growth} therefore gives the objects
$E_d$ and the strict inclusions of their radical tensor ideals for $H$.
We use its notation $\mathcal C_e$ and its bound for $E_e\otimes Y$.

Let $X_d=\operatorname{Ind}_H^G E_d$. Exact induction applied to
the diagonal-cone construction, with the projection formula
\cite[Corollary~4.3.8(4)]{webb2016representations}, gives
$X_e\in\operatorname{Thick}^{\otimes}(X_d)$:
the auxiliary $Q$ is the restriction of $\bigoplus_{j\ge0}R[j]$
with trivial $G$-action.
Mackey's formula \cite[Theorem~5.2.1]{webb2016representations} gives
\[
 \operatorname{Res}_H^G X_d\simeq
        E_d^{\oplus m}\oplus\operatorname{Ind}_1^H W_d,
                      \qquad m\ge1,
\]
for a bounded-above $R$-complex $W_d$ with finitely generated cohomology.
Indeed an intersection of two subgroups of order $p$ is trivial or
the whole subgroup. Normalizer automorphisms send the generator
$u$ of $R/pR$ to a unit multiple, because $R$ is local; hence they
preserve every $P_a$ and $E_d$ up to isomorphism.
On induced objects $u$ acts as zero: $u_A=0$ by projectivity, and
$\operatorname{Ind}_1^H W=A\otimes_R^{\mathbf L}iW$.
The same holds after tensoring with any object.
It follows that
\[
 \operatorname{Thick}^{\otimes}(X_e)
       \subseteq(\operatorname{Res}_H^G)^{-1}\mathcal C_e.
\]
Every power of $\operatorname{Res}_H^G X_d$ contains the corresponding
power of $E_d$ as a retract. By the proof of
Theorem~\ref{thm:abstract-periodic-growth}, these powers do not belong
to $\mathcal C_e$ when $d>e$. Thus the principal radical ideals of
the $X_d$ are strictly nested. 
Forgetting the $H$-action splits $P_a$ as $R\oplus R[2a-1]$.
The projection formula therefore gives
\[
 RG\otimes_R^{\mathbf L}X_d
 \simeq RG\otimes_R^{\mathbf L}\operatorname{Res}_1^G X_d,
 \qquad RG\in\operatorname{Thick}^{\otimes}(X_d),
\]
since the underlying complex of $X_d$ has $R$ as a retract.
Also $H^0X_d\cong\operatorname{Ind}_H^G R$, since $H^0P_1=R$
and $P_a$ has no positive cohomology.
The strict radical chain consequently descends to the quotient
by $\operatorname{Thick}^{\otimes}(RG)$, proving its
nonnoetherianity.
\end{proof}

\subsection{The group-algebra criteria}

\begin{theorem}
\label{thm:all-group-noetherian-ring}
\label{cor:tame-group-noetherian}
\label{cor:group-central-relative}
\label{cor:integral-group-growth}
Let $R$ be commutative noetherian and $G$ finite, with diagonal
 derived tensor over $R$ on $\mathbf D^-_{\mathrm{fg}}(RG)$.
Restriction to $R$ gives a spectral embedding with retraction
\[
 \operatorname{Spc}(\mathbf D^-_{\mathrm{fg}}(R))
 \lhook\joinrel\longrightarrow
 \operatorname{Spc}(\mathbf D^-_{\mathrm{fg}}(RG)),
\]
which is a homeomorphism when $|G|\in R^\times$. Moreover,
\[
 \operatorname{Spc}(\mathbf D^-_{\mathrm{fg}}(RG))\text{ is noetherian}
 \quad\Longleftrightarrow\quad R\text{ is artinian and }|G|\in R^\times.
\]
In that case the spectrum is the finite discrete space $\operatorname{Spec}R$.
For a field $k$, both $\mathbf D^-(kG\text{-}\mathrm{mod})$ and
$\mathbf D^+(kG\text{-}\mathrm{mod})$ have noetherian spectrum exactly
when $|G|\in k^\times$, and their spectra are then singletons.
For every multiplicative $\Phi\subseteq R$,
\[
 \operatorname{Spc}\!\left(
 \left(\frac{\mathbf D^-_{\mathrm{fg}}(RG)}
 {\operatorname{Thick}^{\otimes}(RG)}\right)_\Phi\right)
 \text{ is noetherian}
 \quad\Longleftrightarrow\quad |G|\in(\Phi^{-1}R)^\times.
\]
The algebraic quotient is zero in the invertible case.
\end{theorem}
\begin{proof}
Restriction $F$ and the trivial-action functor $i$ are exact symmetric
monoidal and satisfy $Fi=\operatorname{id}$. Thus
$\operatorname{Spc}(i)\operatorname{Spc}(F)=\operatorname{id}$.
Suppose that $n=|G|$ is invertible in $R$. The norm divided by $n$
splits $R$ from $RG$. For every complex $X$ of $RG$-modules, the maps
\[
 \epsilon_X:RG\otimes_R iF(X)\longrightarrow X,
 \quad g\otimes x\longmapsto gx,
 \qquad
 \sigma_X(x)=\frac1n\sum_{g\in G}g\otimes g^{-1}x
\]
are $G$-linear chain maps and satisfy $\epsilon_X\sigma_X=\operatorname{id}$.
The isomorphism
$RG\otimes_R^{\mathbf L}X\simeq RG\otimes_R^{\mathbf L}iF(X)$
is given degreewise by $g\otimes x\mapsto g\otimes g^{-1}x$.
Since $RG$ is $R$-free, these formulas compute the derived tensor
products and give
\[
 X\mid RG\otimes_R^{\mathbf L}iF(X),\qquad
 iF(X)\mid RG\otimes_R^{\mathbf L}iF(X)
              \simeq RG\otimes_R^{\mathbf L}X.
\]
Consequently
\[
 \operatorname{Thick}^{\otimes}(X)
   =\operatorname{Thick}^{\otimes}(iF(X)),\qquad
 F^{-1}i^{-1}(\mathcal P)=\mathcal P
 \quad(\mathcal P\in\operatorname{Spc}(\mathbf D^-_{\mathrm{fg}}(RG))).
\]
Together with $Fi=\operatorname{id}$, this shows that
$\operatorname{Spc}(F)$ and $\operatorname{Spc}(i)$ are inverse
homeomorphisms.

The retraction and Corollary~\ref{thm:noetherian-balmer-spectrum}
force $R$ to be artinian if the group spectrum is noetherian.
If $R$ is artinian and $|G|\notin R^\times$, projection to a local factor with residue
characteristic $p\mid|G|$ gives a Verdier localization; its spectrum
is nonnoetherian by Proposition~\ref{prop:local-modular-growth}.
This proves the absolute criterion and identifies the spectrum with
the finite discrete space $\operatorname{Spec}R$ in the noetherian case.
Over a field, contragredient duality is an exact symmetric monoidal
anti-equivalence between the two half-bounded categories: tensor
totalizations have finite diagonals. Hence their spectra are homeomorphic.

The norm and augmentation factor $|G|\operatorname{id}_R$ through
$RG$. Thus $|G|$ annihilates the unit in the quotient by
$\operatorname{Thick}^{\otimes}(RG)$, and algebraic localization makes
the quotient zero when $|G|\in(\Phi^{-1}R)^\times$.
Otherwise choose a prime $\mathfrak p$ disjoint from $\Phi$ with
$|G|\in\mathfrak p$. Geometric localization and thick-quotient
interchange give a Verdier localization to
\[
 \mathbf D^-_{\mathrm{fg}}(R_{\mathfrak p}G)/
       \operatorname{Thick}^{\otimes}(R_{\mathfrak p}G)
\]
\cite[Theorem~3.2]{miyachi1991localization}
\cite[Proposition~2.1]{positselski2021unbounded}.
Proposition~\ref{prop:local-modular-growth} gives a strict radical
chain in this quotient. Its spectrum embeds in the original algebraic
quotient by \cite[Proposition~3.11]{balmer2005spectrum}, so the latter
is nonnoetherian.
\end{proof}

For $R=\mathbb Z$, compose the retraction induced by restriction to
$\mathbb Z$ with the retraction associated to
$\mathbb Z\to\mathbb Z_{(p)}$. This gives
$\operatorname{Spc}(\mathbf D^-_{\mathrm{fg}}(\mathbb Z_{(p)}))$
as a spectral retract for every finite $G$ and prime $p$.
\section{Geometric and arithmetic applications}\label{sec:geometric-arithmetic}

\subsection{Schemes and traits}

For a noetherian scheme $X$, $\mathbf D^-_{\mathrm{Coh}}(X)$
is the full subcategory
of $\mathbf D_{\mathrm{qc}}(X)$ with bounded-above coherent
cohomology. Its standard $t$-structure has heart
$\operatorname{coh}X$.

For an open immersion $j:U\hookrightarrow X$ with complement $Z$, the
coherent localization, obtained from Gabriel's theorem by passing
to noetherian objects
\cite[Chapter~VI, Section~1, Proposition~3]{gabriel1962}
\cite[Tag~01PF]{stacks-project}, is
\[
 (\operatorname{coh}X)/(\operatorname{coh}_Z X)
             \simeq\operatorname{coh}U.
\]
Miyachi's localization theorem
\cite[Theorem~3.2]{miyachi1991localization}, in the formulation of
\cite[Theorem~1.1]{wuzhou2022shortexact}, and
\cite[Tag~0FDB]{stacks-project} therefore give
\[
 \mathbf D^-_{\mathrm{Coh}}(X)/\mathbf D^-_{\mathrm{Coh},Z}(X)
                  \simeq\mathbf D^-_{\mathrm{Coh}}(U).
\]
Thus restriction induces a topological embedding on Balmer spectra
\cite[Proposition~3.11]{balmer2005spectrum}.

\begin{corollary}
\label{cor:scheme-noetherianity-criterion}
\label{cor:scheme-valuation-retraction}
For a nonempty noetherian scheme $X$,
\[
 \operatorname{Spc}\bigl(\mathbf D^-_{\mathrm{Coh}}(X)\bigr)
       \text{ is noetherian}
 \quad\Longleftrightarrow\quad \dim X=0.
\]
If $\dim X>0$, there is a morphism $v:\operatorname{Spec}V\to X$
from a discrete valuation ring such that
\[
 \operatorname{Spc}(Lv^*):
 \operatorname{Spc}\bigl(\mathbf D^-_{\mathrm{fg}}(V)\bigr)
 \hookrightarrow
 \operatorname{Spc}\bigl(\mathbf D^-_{\mathrm{Coh}}(X)\bigr)
\]
is a topological embedding.
\end{corollary}
\begin{proof}
If $\dim X>0$, choose distinct points $y\leadsto x$ and an affine
neighbourhood $U=\operatorname{Spec}A$ of $x$. It contains $y$, so $A$
is not artinian. The construction in
Corollary~\ref{thm:noetherian-balmer-spectrum} gives a DVR $V$ and
$A\to V$ for which
$\operatorname{Spc}(\mathbf D^-_{\mathrm{fg}}(V))$ embeds in
$\operatorname{Spc}(\mathbf D^-_{\mathrm{fg}}(A))$.
Compose this embedding with the embedding of
the restriction functor above.
The DVR spectrum is not noetherian, and every subspace of a noetherian
space is noetherian. Thus the spectrum for $X$ is not noetherian.
If $\dim X=0$, its underlying noetherian space is finite and discrete.
Each singleton is an affine open with artinian coordinate ring.
Hence $X$ is affine with artinian coordinate ring, and the conclusion
follows from Corollary~\ref{thm:noetherian-balmer-spectrum}.
\end{proof}

If $s\in\Gamma(X,\mathcal O_X)$ and
$u:\operatorname{Spec}V\to X$ satisfy $u^*s=\varepsilon\pi^c$ with
$c\ge1$, then the embedding strengthens to a spectral retraction by
Theorem~\ref{thm:valuation-growth-retraction}. Indeed derived
pullback preserves pseudo-coherent objects, tensor products and the
coproducts $\bigoplus_{n\ge0}\mathcal O_X[n]^{d_n}$
\cite[Tags~08E8, 09U7, 09J3 and 08EU]{stacks-project}. Apply the theorem to the
objects $\bigoplus_{n\ge0}\operatorname{Cone}(s^{f(n)})[n]$.

\begin{theorem}\label{thm:central-relative-scheme}
Let $X$ be a noetherian scheme, let
$\Phi\subseteq\Gamma(X,\mathcal O_X)$ be multiplicative, and suppose
that $\mathcal C=\operatorname{Thick}(\mathcal E)$ is a tensor ideal
of $\mathbf D^-_{\mathrm{Coh}}(X)$. Write
\[
 W=\bigcup_{E\in\mathcal E,i}\operatorname{Supp}H^iE,
 \qquad
 U=(X\setminus W)\cap\bigcap_{\phi\in\Phi}D(\phi).
\]
Suppose that $\mathbf D^{\mathrm{perf}}_W(X)\subseteq\mathcal C$.
For each nonclosed $y\in U$, suppose there exists a morphism
$v_y:\operatorname{Spec}V_y\to X$ from a DVR, whose generic image
is $y$ and whose closed image differs from $y$, such that
\begin{equation}\label{eq:relative-trait-uniformity}
 \forall E\in\mathcal E\quad\exists a(E,y)\ge1:\qquad
 \pi_y^{a(E,y)}\operatorname{id}_{\mathbf L v_y^*E}=0.
\end{equation}
Then
\begin{equation}\label{eq:central-relative-noetherianity}
 \operatorname{Spc}\!\left(
     (\mathbf D^-_{\mathrm{Coh}}(X)/\mathcal C)_\Phi\right)
       \text{ is noetherian}
 \quad\Longleftrightarrow\quad
       \text{every point of }U\text{ is closed in }X.
\end{equation}
In this case $U$ is finite and clopen and
\[
 (\mathbf D^-_{\mathrm{Coh}}(X)/\mathcal C)_\Phi
       \simeq\mathbf D^-_{\mathrm{Coh}}(U).
\]
For $X=\operatorname{Spec}R$, the perfect-object inclusion is
automatic.
\end{theorem}
\begin{proof}
Write $\mathcal T=\mathbf D^-_{\mathrm{Coh}}(X)$ and
\[
 W'=W\cup\bigcup_{\phi\in\Phi}Z(\phi),\qquad
 \mathcal D=\operatorname{Thick}(\mathcal C\cup\mathcal S_\Phi).
\]
Algebraic central localization commutes with the thick quotient:
\begin{equation}\label{eq:central-relative-commutation}
 (\mathcal T/\mathcal C)_\Phi
 \simeq\mathcal T/\mathcal D
 \simeq\mathcal T_\Phi/
         \operatorname{Thick}(q_\Phi\mathcal C).
\end{equation}
The description of $\mathcal S_\Phi$ by cones of scalars and Thomason's
perfect-support classification give
\begin{align}
 \mathcal S_\Phi
   &=\operatorname{Thick}^{\otimes}
        \{\operatorname{Cone}(\phi:\mathcal O_X\to\mathcal O_X):
                                                       \phi\in\Phi\},
       \label{eq:relative-central-cones}\\
 \mathbf D^{\mathrm{perf}}_{W'}(X)
   &\subseteq\mathcal D
     \subseteq\mathbf D^-_{\mathrm{Coh},W'}(X).
       \label{eq:relative-central-sandwich}
\end{align}

Suppose that $y\in U$ is not closed. Choose the trait $v_y$ in the
hypotheses and an affine open $j:X_0\hookrightarrow X$ containing its
closed image. The trait factors through $X_0$. Restriction is a
Verdier localization, so quotient interchange gives a Verdier localization
\[
 (\mathcal T/\mathcal C)_\Phi\longrightarrow
 \left(\frac{\mathbf D^-_{\mathrm{Coh}}(X_0)}
 {\mathcal C_0}\right)_{j^*\Phi},\qquad
 \mathcal C_0=\operatorname{Thick}(j^*\mathcal E).
\]
Since $j^*$ is essentially surjective and $\mathcal C$ is a tensor
ideal, $\mathcal C_0$ is a tensor ideal. Derived pullback
$F=\mathbf L v_y^*$ on $X_0$ preserves pseudo-coherent objects and
coproducts of shifts of the unit with finite multiplicity in each
degree, as above. Write $\mathfrak p\subsetneq
\mathfrak q$ for the generic and closed images of the trait in $X_0$.
Choose $s\in\mathfrak q\setminus\mathfrak p$; then
\[
 F(s)=\varepsilon\pi_y^c\quad(c\ge1),\qquad
 F(j^*\phi)\ne0\quad(\phi\in\Phi),\qquad
 \pi_y^{a(E,y)}\operatorname{id}_{F(j^*E)}=0\quad(E\in\mathcal E).
\]
The second condition follows from $y\in D(\phi)$, and the third is
\eqref{eq:relative-trait-uniformity}.
Corollary~\ref{thm:central-relative-growth} therefore gives strictly
nested radical thick tensor ideals, and a strictly descending chain
of closed supports, in the
target quotient. Its spectrum embeds in
$\operatorname{Spc}((\mathcal T/\mathcal C)_\Phi)$ by
\cite[Proposition~3.11]{balmer2005spectrum}. Since every subspace of a
noetherian space is noetherian, the latter spectrum is not noetherian.

Conversely, $U=X\setminus W'$ is generalization closed. If every
point of $U$ is closed in $X$, these points are singleton irreducible
components disjoint from the other components. Noetherianity makes
$U$ finite and clopen. The perfect idempotent $\mathcal O_{W'}$
belongs to $\mathcal D$ by
\eqref{eq:relative-central-sandwich}, so
\[
 \mathcal D=\mathbf D^-_{\mathrm{Coh},W'}(X),\qquad
 \mathcal T/\mathcal D\simeq\mathbf D^-_{\mathrm{Coh}}(U).
\]
The scheme $U$ is affine and artinian. Its spectrum is noetherian
by the zero-dimensional case.

Finally, if $X=\operatorname{Spec}R$, Matsui--Takahashi
\cite[Proposition~2.9 and Corollary~2.11]{Matsui_2017} show that
every tensor ideal with support $W$ contains all perfect complexes
supported on $W$.
\end{proof}

The exponent in~\eqref{eq:relative-trait-uniformity} may depend on
both generator and trait. Taking $\mathcal E=\varnothing$ gives the
algebraic central-localization criterion; taking $\Phi=\{1\}$ gives
the relative criterion without scalar localization.
For an ordinary thick subcategory that is not a tensor ideal,
\eqref{eq:central-relative-commutation} remains a triangulated
equivalence, but the indicated Balmer spectra are not defined.

\begin{example}\label{ex:unbounded-relative-central}
Let $R=\overline{\mathbb Q}[t]$, enumerate
$\overline{\mathbb Q}=\{a_1,a_2,\ldots\}$, and let
\[
 E=\bigoplus_{n\ge1}R/(t-a_n)^n[n],\qquad
 \mathcal C=\operatorname{Thick}^{\otimes}(E),\qquad
 \Phi=R\setminus\{0\}.
\]
Equivalently $\mathcal C=\operatorname{Thick}(\mathcal E)$ for
$\mathcal E=\{E\otimes_R^{\mathbf L}Y:Y\in\mathbf D^-_{\mathrm{fg}}(R)\}$.
Here $W$ is the set of closed points and $U=\{(0)\}$.
On every trait with generic image $(0)$, at most one summand of
$E$ survives, so \eqref{eq:relative-trait-uniformity} holds.
Every cone of a scalar is perfect with support in $W$, whence
$\mathcal S_\Phi\subseteq\mathcal C$ and
\[
 (\mathbf D^-_{\mathrm{fg}}(R)/\mathcal C)_\Phi
       \simeq\mathbf D^-_{\mathrm{fg}}(R)/\mathcal C.
\]
Its spectrum is not noetherian by
\eqref{eq:central-relative-noetherianity}. By contrast, geometric
localization at $\Phi$ gives
$\mathbf D^-_{\mathrm{fg}}(\overline{\mathbb Q}(t))$, whose spectrum
is a point. Moreover
\[
 \operatorname{Thick}^{\otimes}(\mathbf D^{\mathrm{perf}}_W(R))
       \subsetneq\mathcal C
       \subsetneq\mathbf D^-_{\mathrm{fg}}(R)_W.
\]
The first strictness follows from the infinite support of $E$;
the second follows by localizing
$\bigoplus_{m\ge1}R/(t-a_1)^m[m]$ at $(t-a_1)$, where it has no
uniform annihilator while every localized object of $\mathcal C$
does.
\end{example}

\subsection{Arithmetic spectra}
For $R=S^{-1}\mathbb Z\ne0$, write $\mathcal K_R$ for the $S$-local
spectra bounded below with degreewise finitely generated homotopy
over $R$, with smash product and the Postnikov $t$-structure.

\begin{corollary}
\label{cor:arithmetic-spectra-noetherian}
\label{cor:arithmetic-central-relative}
\label{cor:t-structured-finite-type-spectra}
The category $\mathcal K_R$ is essentially small and tensor triangulated,
and
\[
 \operatorname{Spc}(\mathcal K_R)\text{ is noetherian}
 \quad\Longleftrightarrow\quad R\text{ is artinian}
 \quad\Longleftrightarrow\quad R=\mathbb Q.
\]
Over $\mathbb Q$ its spectrum is a singleton. Every prime $p$ not
invertible in $R$ gives a spectral retract
\[
 \operatorname{Spc}(\mathbf D^-_{\mathrm{fg}}(\mathbb Z_{(p)}))
 \lhook\joinrel\longrightarrow\operatorname{Spc}(\mathcal K_R).
\]
Let $\mathcal C=\operatorname{Thick}(\mathcal E)$ be a tensor ideal
of $\mathcal K_{\mathbb Z_{(p)}}$ such that for every $E\in\mathcal E$
there is $a(E)\ge1$ with $p^{a(E)}\operatorname{id}_E=0$. Then
\[
 (\mathcal K_{\mathbb Z_{(p)}}/\mathcal C)_{\{p^a:a\ge0\}}
 \simeq\mathcal K_{\mathbb Z_{(p)}}/\mathcal S_p,
 \qquad
 \operatorname{Spc}(\mathcal K_{\mathbb Z_{(p)}}/\mathcal S_p)
 \text{ is not noetherian}.
\]
Geometric rationalization has target $\mathcal K_{\mathbb Q}$.
\end{corollary}
\begin{proof}
Write $A=\mathbb S[S^{-1}]$. Serre's finiteness theorem
\cite{serre1953homotopie} makes $A$ coherent; the almost-perfect criteria and finite-free resolutions in
\cite[Propositions~7.2.4.11(5), 7.2.4.17--7.2.4.18]{lurie2017ha}
give the asserted category properties. For every noninvertible $p$,
extension to $H\mathbb Z_{(p)}$ preserves almost-perfect modules and
provides the exact symmetric monoidal functor
\[
 F_p(E)=H\mathbb Z_{(p)}\wedge_A E:
 \mathcal K_R\longrightarrow\mathbf D^-_{\mathrm{fg}}(\mathbb Z_{(p)})
\]
\cite[Theorem~7.1.2.13]{lurie2017ha}.
It preserves the coproducts $P_{\mathbf d}$ and takes $p$ to a DVR
uniformizer. Theorem~\ref{thm:valuation-growth-retraction} gives the
spectral retract and nonnoetherianity.
If $R=\mathbb Q$, Serre's theorem \cite{serre1953homotopie} identifies
$A$ with $H\mathbb Q$. Thus
$\mathcal K_R\simeq\mathbf D^-_{\mathrm{fg}}(\mathbb Q)$, where every
nonzero object has a shifted unit as a retract. Its spectrum is a point.
A nonzero localization of $\mathbb Z$ is artinian exactly when it is
$\mathbb Q$.

For the relative assertion, $\mathcal C\subseteq\mathcal S_p$.
Proposition~\ref{prop:minimal-central-localization} gives the
equivalence, and
Corollary~\ref{thm:central-relative-growth}, applied to $F_p$ and $p$,
gives a strict chain after quotienting by $\mathcal S_p$.
The rational target is the case just proved.
\end{proof}

For $\Phi=\{p^a:a\ge0\}$ in $\mathcal K_{\mathbb Z}$,
Proposition~\ref{thm:abstract-central-t-structure} gives
$\operatorname{PC}^-(\mathcal S_\Phi)=\mathcal N_\Phi$.
The object
\[
 E=\bigvee_{n\ge0}\operatorname{Cone}
                (p^{n+1}:\mathbb S\to\mathbb S)[n]
\]
has degreewise $p$-power-torsion homotopy but no uniform annihilator:
its $n$th summand has $\pi_n=\mathbb Z/p^{n+1}$.
Thus the central quotient has heart
$\mathbb Z[1/p]\text{-}\mathrm{mod}$ and nonzero objects with zero
cohomology in every degree.

Nonnoetherianity also follows from classical chromatic theory.
The fully faithful inclusion of finite $S$-local spectra into
$\mathcal K_R$ detects tensor-nilpotence of morphisms, so its map on
spectra is surjective \cite[Theorem~1.3]{balmer2018surjectivity}.
For $p$ not invertible in $R$, the chromatic chain gives a nonnoetherian
subspace of the target \cite[Theorem~9.1 and Corollary~9.5]{balmer2010spectra}.
The corollary also gives a DVR spectral retraction and a strictly
nested family of radical thick tensor ideals after quotienting by
$\mathcal S_p$.
\section{Cohomological information and localization}
\label{sec:cohomological-information}
\label{sec:hom-modules}\label{sec:postnikov-tails}\label{sec:balmer}

The radical thick tensor ideals constructed above have the same
cohomological support.
We describe the information retained by their cohomology, and the two
kernels that arise when localization is compared with thick generation.

\subsection{The heart and Postnikov closure}
Let $\mathcal T$ be an essentially small triangulated category with a
$t$-structure, heart $\mathcal A$, and set
$\mathcal T^- =\bigcup_b\mathcal T^{\le b}$,
$\mathcal T^+ =\bigcup_a\mathcal T^{\ge a}$ and
$\mathcal T^b=\mathcal T^-\cap\mathcal T^+$.
A thick subcategory of $\mathcal T^-$ is \emph{Postnikov-complete}
if it is truncation-closed and contains $X$ whenever it contains
$\tau^{\ge n}X$ for every $n$; in $\mathcal T^+$ use $\tau^{\le n}$.
Write $\operatorname{PC}^\pm$ for this closure.
A full additive subcategory of $\mathcal A$ is \emph{wide} if it is
closed under kernels, cokernels and extensions.

\begin{proposition}
\label{thm:abstract-postnikov-heart}\label{def:postnikov-completeness}
\label{cor:classifying-heart-support}
Intersection with $\mathcal A$ identifies Postnikov-complete thick
subcategories of $\mathcal T^\pm$ with wide subcategories of $\mathcal A$.
In particular,
\begin{equation}\label{eq:abstract-postnikov-heart}
 \operatorname{PC}^\pm(\mathcal E)
 =\{X\in\mathcal T^\pm:H^iX\in
    \operatorname{Wide}_{\mathcal A}\{H^jE:E\in\mathcal E,\ j\in\mathbb Z\}
                                      \ (\forall i)\}.
\end{equation}
Every truncation-closed thick $\mathcal Z\subseteq\mathcal T^\pm$ satisfies
\begin{equation}\label{eq:abstract-bounded-heart}
 \mathcal Z\cap\mathcal T^b
   =\{X\in\mathcal T^b:H^iX\in\mathcal Z\cap\mathcal A\ (\forall i)\}.
\end{equation}
\end{proposition}
\begin{proof}
For a wide $\mathcal W\subseteq\mathcal A$, the sequences
\[
 0\longrightarrow\operatorname{coker}H^i(f)
 \longrightarrow H^i\operatorname{Cone}(f)
 \longrightarrow\ker H^{i+1}(f)\longrightarrow0
\]
show that $\{X:H^iX\in\mathcal W\ (\forall i)\}$ is thick.
It is truncation-closed and Postnikov-complete.
The bounded correspondence \cite[Theorem~2.5]{ZhangCai2017} gives
\eqref{eq:abstract-bounded-heart} and makes $\mathcal Z\cap\mathcal A$ wide.
If $\mathcal Z$ is Postnikov-complete, then
\[
 X\in\mathcal Z
 \ \Longleftrightarrow\ \tau^{\ge n}X\in\mathcal Z\cap\mathcal T^b\ (\forall n)
 \ \Longleftrightarrow\ H^iX\in\mathcal Z\cap\mathcal A\ (\forall i)
\]
in $\mathcal T^-$; use $\tau^{\le n}$ in $\mathcal T^+$.
Taking the least wide subcategory containing the $H^jE$ proves the formula.
\end{proof}

If a support $\sigma$ classifies wide subcategories of $\mathcal A$ by
specialization-closed subsets, extend it by $\sigma(X)=\bigcup_i\sigma(H^iX)$.
Then
\begin{equation}\label{eq:abstract-classifying-heart-closure}
 \operatorname{PC}^\pm(\mathcal E)=\mathcal T^\pm_{\sigma(\mathcal E)}.
\end{equation}
For truncation-closed thick $\mathcal Z$, also
$\mathcal Z\cap\mathcal T^b=\mathcal T^b_{\sigma(\mathcal Z)}$.
If $\mathcal P\subseteq\mathcal T^b$ realizes every such support, then
$\mathcal Z\cap\mathcal P=\mathcal P_{\sigma(\mathcal Z)}$ and
$\sigma(\mathcal Z\cap\mathcal P)=\sigma(\mathcal Z)$.
These conclusions require a classification in the heart, as the following
example shows.

\begin{example}
\label{ex:projective-wide-postnikov}
Let $X=\mathbb P^1_k$. In either half-bounded coherent category,
Proposition~\ref{thm:abstract-postnikov-heart} gives
\[
 \operatorname{PC}^{\pm}
       (\operatorname{Thick}(\mathcal O_X))
   =\{E\mid H^iE\in\operatorname{add}(\mathcal O_X)
                                               \ (\forall i)\}
   \subsetneq\mathbf D^{\pm}_{\mathrm{Coh}}(X).
\]
Indeed $\operatorname{Hom}(\mathcal O_X,\mathcal O_X)=k$ and
$\operatorname{Ext}^1(\mathcal O_X,\mathcal O_X)=0$
\cite[Tag~01XT]{stacks-project}; matrices of constants have
kernels and cokernels in $\operatorname{add}(\mathcal O_X)$,
so this is a wide subcategory. It is not Serre, since
$\mathcal O_X(-1)\hookrightarrow\mathcal O_X$ and
$\mathcal O_X(-1)\notin\operatorname{add}(\mathcal O_X)$.
Both sides of the strict inclusion have full ordinary support.
These are ordinary thick subcategories: the tensor ideal
generated by $\mathcal O_X$ is the whole category.
\end{example}

\subsection{Central localization and cohomology}
Let a commutative ring $C$ act through the degree-zero center of a triangulated
category $\mathcal T$ with heart $\mathcal A$, and let
$\Phi\subseteq C$ be multiplicative. In either half-bounded category, write
\[
\begin{aligned}
 \mathcal S_\Phi^\pm&=\{X:\phi\operatorname{id}_X=0
                                        \text{ for some }\phi\in\Phi\},
 &\mathcal A_\Phi&=\mathcal S_\Phi^\pm\cap\mathcal A,\\
 \mathcal N_\Phi^\pm&=\{X:H^iX\in\mathcal A_\Phi\ (\forall i)\},
 &\mathcal Q_\Phi^\pm&=\mathcal T^\pm/\mathcal S_\Phi^\pm.
\end{aligned}
\]

\begin{proposition}\label{prop:central-comparison}
\label{thm:abstract-central-t-structure}
\label{cor:abstract-geometric-comparison}
\label{lem:early-negative-localization}
The images of the truncation classes define a $t$-structure on
$\mathcal Q_\Phi^\pm$. If $q$ is the quotient and
$\pi:\mathcal A\to\mathcal A/\mathcal A_\Phi$, then
\begin{gather}
 (\mathcal Q_\Phi^\pm)^\heartsuit\simeq\mathcal A/\mathcal A_\Phi,
 \qquad H^i(qX)=\pi H^iX,\qquad
 \operatorname{PC}^\pm(\mathcal S_\Phi^\pm)=\mathcal N_\Phi^\pm,
       \label{eq:abstract-algebraic-geometric-completion}\\
 \mathcal D_\Phi^-:=\bigcap_n(\mathcal Q_\Phi^-)^{\le n}
       \simeq\mathcal N_\Phi^-/\mathcal S_\Phi^-,\qquad
 \mathcal D_\Phi^+:=\bigcap_n(\mathcal Q_\Phi^+)^{\ge n}
       \simeq\mathcal N_\Phi^+/\mathcal S_\Phi^+.
       \label{eq:central-comparison-degenerate-kernel}
\end{gather}
The opposite intersections are zero, and
\[
 \mathcal D_\Phi^-=\{qX:q\tau^{\ge n}X=0\ (\forall n)\},
 \qquad
 \mathcal D_\Phi^+=\{qX:q\tau^{\le n}X=0\ (\forall n)\}.
\]
Moreover,
\begin{gather*}
 \mathcal Q_\Phi^\pm\text{ is nondegenerate}
 \quad\Longleftrightarrow\quad
 \mathcal S_\Phi^\pm=\mathcal N_\Phi^\pm\\
 \Longleftrightarrow\quad
 \text{bounded standard truncations detect zero objects}.
\end{gather*}

Suppose $L:\mathcal T^\pm\to\mathcal V^\pm$ is a $t$-exact Verdier
localization, target cohomology detects zero objects, and
$\ker L^\heartsuit=\mathcal A_\Phi$. Then
\[
 \ker L=\mathcal N_\Phi^\pm,\qquad
 \overline L:\mathcal Q_\Phi^\pm\longrightarrow\mathcal V^\pm,
 \qquad\ker\overline L\simeq\mathcal N_\Phi^\pm/\mathcal S_\Phi^\pm.
\]
Thus $\overline L$ is an equivalence if and only if
$\mathcal S_\Phi^\pm=\mathcal N_\Phi^\pm$.
\end{proposition}
\begin{proof}
Scalar annihilators are preserved by truncation. The annihilator
product rule makes $\mathcal S_\Phi^\pm$ thick and
$\mathcal A_\Phi$ Serre. Since
\[
 \mathcal A_\Phi\subseteq\mathcal S_\Phi^\pm,
 \qquad H^i(\mathcal S_\Phi^\pm)\subseteq\mathcal A_\Phi,
\]
Proposition~\ref{thm:abstract-postnikov-heart} gives the Postnikov equality.
Scalar fractions give $t$-orthogonality in the quotient and identify
its heart with the Gabriel quotient; truncation triangles give
$H^iq=\pi H^i$.

For a bounded $B\in\mathcal T^{[a,b]}$,
\[
 \prod_{i=a}^b\operatorname{Ann}_C H^iB
 \subseteq\operatorname{Ann}_C B
 \subseteq\bigcap_{i=a}^b\operatorname{Ann}_C H^iB.
\]
Hence $qB=0$ precisely when all its cohomology lies in
$\mathcal A_\Phi$. Applying this to the bounded truncations gives
the formulas for $\mathcal D_\Phi^\pm$; half-boundedness and
orthogonality give the opposite zero intersections.

For $L$ as stated, $LX=0$ if and only if
$L^\heartsuit H^iX=0$ for every $i$. This identifies its kernel
and Verdier transitivity gives
$\ker\overline L\simeq\mathcal N_\Phi^\pm/\mathcal S_\Phi^\pm$
\cite[Tag~05RJ]{stacks-project}. The kernel formulas prove the
equivalence criteria.
\end{proof}

For a commutative noetherian ring $R$, extension of scalars gives the
three Verdier localizations
\[
\begin{array}{c|c}
 \text{source}&\text{target}\\ \hline
 \mathbf D^-_{\mathrm{fg}}(R)&\mathbf D^-_{\mathrm{fg}}(\Phi^{-1}R)\\
 \mathbf D^+_{\mathrm{fg}}(R)&\mathbf D^+_{\mathrm{fg}}(\Phi^{-1}R)\\
 \mathbf D^+(R\text{-}\mathrm{mod})
    &\mathbf D^+(\Phi^{-1}R\text{-}\mathrm{mod}).
\end{array}
\]
The first is \cite[Proposition~6.18]{sanders2026balmer}; the two positive
realizations are proved separately in
Proposition~\ref{prop:geometric-localization-plus}.
Their heart functors are
$(R\text{-}\mathrm{mod})\to(\Phi^{-1}R)\text{-}\mathrm{mod}$
by Corollary~\ref{cor:finite-module-central-localization}.
With $\mathcal S_\Phi$ defined by a uniform identity annihilator and
$\mathcal N_\Phi$ by degreewise torsion, the comparison kernels are
\begin{align}
 \ker\overline L^-_\Phi&\simeq\mathcal N^-_\Phi/\mathcal S^-_\Phi,
       \label{eq:minus-degenerate-tail}\\
 \ker\overline L^+_\Phi&\simeq\mathcal N^+_\Phi/\mathcal S^+_\Phi,
       \label{eq:plus-degenerate-tail}\\
 \ker\overline L^{+,\mathrm{ft}}_\Phi
   &\simeq\mathcal N^{+,\mathrm{ft}}_\Phi/\mathcal S^{+,\mathrm{ft}}_\Phi.
       \label{eq:finite-term-comparison-defect}
\end{align}
Thus Proposition~\ref{prop:central-comparison} applies to all three rows.

\subsection{Derived Hom and perfect generation}
Throughout the rest of this section $R$ is commutative noetherian.
Here $\mathbf D^+_{\mathrm{fg}}(R)$ is the full subcategory of
$\mathbf D(R\text{-}\mathrm{Mod})$ with bounded-below finite cohomology;
$\mathbf D^+(R\text{-}\mathrm{mod})$ is the derived category of bounded-below
complexes of finite modules. We use these two realizations separately.
All constructions in $\mathbf D^+(R\text{-}\mathrm{mod})$ are internal.
Its canonical functor to $\mathbf D^+_{\mathrm{fg}}(R)$ is denoted by $j_R$.
Actions here are coherent unital biexact actions, and thick generation
uses finite triangles and retracts
\cite[Definition~3.2(1)--(3)]{Stevenson}.

\begin{lemma}\label{prop:derived-hom-module}\label{lem:hom-generation-normal-form}
The categories $\mathbf D^+_{\mathrm{fg}}(R)$ and
$\mathbf D^+(R\text{-}\mathrm{mod})$ are triangulated modules over
$(\mathbf D^-_{\mathrm{fg}}(R))^{\mathrm{op}}$ under
$A\odot X=\mathbf R\operatorname{Hom}_R(A,X)$.
The functor $j_R$ is exact, conservative, and compatible with the actions
and standard truncations. In either positive category,
\begin{equation}\label{eq:hom-generation-normal-form}
 \operatorname{Thick}^{\operatorname{Hom}}(\mathcal E)
   =\operatorname{Thick}\{A\odot E:A\in\mathbf D^-_{\mathrm{fg}}(R),E\in\mathcal E\},
 \qquad \operatorname{Supp}(A\odot X)\subseteq\operatorname{Supp}X.
\end{equation}
\end{lemma}
\begin{proof}
Appendix~\ref{app:internal-hom} constructs the action and verifies
its coherence. Currying and exactness give
\eqref{eq:hom-generation-normal-form}. The sums in
\eqref{eq:finite-term-hom-diagonal} are finite, so localization
commutes with them and gives the support inclusion.
\end{proof}

A thick submodule for this action is called a \emph{thick Hom-submodule}.
Write
\[
 \mathcal H_R=\operatorname{Thick}^{\operatorname{Hom}}(R),\qquad
 \mathcal F^+_W=\operatorname{Thick}^{\operatorname{Hom}}
                         (\mathbf D^{\mathrm{perf}}(R)_W).
\]
Use $\mathcal H_R^{+,\mathrm{ft}}$ and $\mathcal F^{+,\mathrm{ft}}_W$ in
$\mathbf D^+(R\text{-}\mathrm{mod})$, and compute identity annihilators
in that category. Both positive categories have the same bounded and
perfect subcategories, through the usual bounded equivalence.

\begin{proposition}\label{thm:perfect-hom-classification}
\label{def:perfect-hom-generation}\label{cor:closed-perfect-hom-submodule}
\label{lem:derived-scalar-annihilators}
Support classifies thick Hom-submodules generated by perfect complexes:
\begin{equation}\label{eq:perfect-hom-intrinsic}
 \mathcal F^+_W
   =\{X\in\mathcal H_R:V(\operatorname{Ann}_{\mathbf D(R)}X)\subseteq W\},
 \qquad
 \mathcal F^+_W\cap\mathbf D^{\mathrm{perf}}(R)=\mathbf D^{\mathrm{perf}}(R)_W.
\end{equation}
For $\mathcal E\subseteq\mathbf D^{\mathrm{perf}}(R)$,
$\operatorname{Thick}^{\operatorname{Hom}}(\mathcal E)=
\mathcal F^+_{\operatorname{Supp}\mathcal E}$.
In particular,
\[
 \mathcal F^+_{V(I)}
   =\{X\in\mathcal H_R:I^n\operatorname{id}_X=0\text{ for some }n\ge1\}.
\]
The same statements hold internally for
$\mathcal F_W^{+,\mathrm{ft}}\subseteq\mathbf D^+(R\text{-}\mathrm{mod})$.
\end{proposition}
\begin{proof}
For a triangle $X\to Y\to Z$, identity annihilators satisfy
$\operatorname{Ann}X\operatorname{Ann}Z\subseteq\operatorname{Ann}Y$;
for bounded finite cohomology,
$V(\operatorname{Ann}X)=\operatorname{Supp}X$
\cite[Proposition~2.3(2) and Lemma~2.5(1)]{Matsui_2017}.
Thus the right-hand side of \eqref{eq:perfect-hom-intrinsic} is a thick
Hom-submodule containing the required perfect objects.
Conversely, let $J=\operatorname{Ann}X=(a_1,\ldots,a_m)$.
The scalar triangles split, giving
\[
 K(J)\otimes_R^{\mathbf L}X
   \simeq\bigoplus_{i=0}^m X[i]^{\binom mi},\qquad
 X\mid K(J)\odot X.
\]
The second assertion uses perfect duality and a shift.
By \eqref{eq:hom-generation-normal-form},
\[
 K(J)\odot\mathcal H_R\subseteq
 \operatorname{Thick}^{\operatorname{Hom}}(K(J)^\vee).
\]
If $V(J)\subseteq W$, this contains $X$ and lies in $\mathcal F^+_W$.
The classification of perfect subcategories
\cite[Theorem~3.15]{Thomason} gives the remaining assertions, since
$V(J)\subseteq V(I)$ is equivalent to $I^n\subseteq J$ for some $n$.
All splittings and annihilators can be read internally for finite terms.
\end{proof}

\begin{corollary}\label{prop:artinian-hom-submodules}
If $R$ is artinian, then $\mathcal H_R=\mathbf D^+_{\mathrm{fg}}(R)$,
and every thick Hom-submodule is $\mathcal F^+_W=
\mathbf D^+_{\mathrm{fg}}(R)_W$ for a unique $W\subseteq\operatorname{Spec}R$.
The same holds in $\mathbf D^+(R\text{-}\mathrm{mod})$.
\end{corollary}
\begin{proof}
Suppose first that $(R,\mathfrak m,k)$ is local and $0\ne Y$.
For $a=\inf\{i:H^iY\ne0\}$,
\[
 H^a\mathbf R\operatorname{Hom}_R(k,Y)
   =\operatorname{Hom}_R(k,H^aY)\ne0,
 \qquad k[-a]\mid\mathbf R\operatorname{Hom}_R(k,Y).
\]
The retract follows by splitting the complex of $k$-vector spaces.
Every bounded-below complex of degreewise finite-dimensional
$k$-vector spaces is
\[
 \bigoplus_{i\ge b}k^{d_i}[-i]
  \simeq\mathbf R\operatorname{Hom}_R
       \left(\bigoplus_{i\ge b}R^{d_i}[i],k\right).
\]
A bounded-below representative $Q$ of finite modules exists by
Remark~\ref{rem:finite-term-comparison}. The finite filtration by
$\mathfrak m^jQ$, with $\mathfrak m^t=0$, now puts every such object
in $\operatorname{Thick}^{\operatorname{Hom}}(k)$.
For $R=\prod R_i$, the functors $e_iR\odot(-)=e_i(-)$ split every
Hom-submodule into the local factors; each factor is zero or the whole
category. Finally $e_WR$ is a perfect generator of the support part.
The argument is internal in the finite-term category as well.
\end{proof}

\subsection{Support and the two localization kernels}
For specialization-closed $W\subseteq\operatorname{Spec}R$, write
$\mathcal G^-_W=\operatorname{Thick}^{\otimes}(\mathbf D^{\mathrm{perf}}(R)_W)$.

\begin{proposition}\label{thm:postnikov-support}
\label{prop:truncation-perfect-detection}\label{thm:uniform-generation-completion}
\label{cor:three-support-classifications}
For every family $\mathcal E$ in either half-bounded ring category,
\begin{equation}\label{eq:postnikov-support-closure}
 \operatorname{PC}^\pm(\mathcal E)
   =\mathbf D^\pm_{\mathrm{fg}}(R)_{\operatorname{Supp}\mathcal E}.
\end{equation}
In the internal finite-term category,
\begin{equation}\label{eq:finite-term-postnikov-support}
 \operatorname{PC}^{+,\mathrm{ft}}(\mathcal E)
   =(\mathbf D^+(R\text{-}\mathrm{mod}))_{\operatorname{Supp}\mathcal E}.
\end{equation}
In particular,
\begin{equation}\label{eq:uniform-generation-completion}
 \operatorname{PC}^-(\mathcal G^-_W)=\mathbf D^-_{\mathrm{fg}}(R)_W,
 \qquad \operatorname{PC}^+(\mathcal F^+_W)=\mathbf D^+_{\mathrm{fg}}(R)_W.
\end{equation}
For a thick tensor ideal $\mathcal J$ one has
\begin{equation}\label{eq:postnikov-tame}
 \operatorname{PC}^-(\mathcal J)=\mathcal J^{\mathrm{tame}}.
\end{equation}
\end{proposition}
\begin{proof}
The three hearts are $R\text{-}\mathrm{mod}$. Its wide subcategories are
the support Serre subcategories
\cite[Corollary~3.2 and Theorem~3.6]{Takahashi}; apply
Proposition~\ref{thm:abstract-postnikov-heart}.
The perfect generators realize $W$ by Koszul complexes.
The last formula is \cite[Lemma~5.8(1)]{Matsui_2017}.
\end{proof}

For a multiplicative $\Phi\subseteq R$, let
$W_\Phi=\bigcup_{\phi\in\Phi}V(\phi)$.
In each category define $\mathcal S_\Phi$ by
$\phi\operatorname{id}_X=0$ for some $\phi\in\Phi$, and
$\mathcal N_\Phi$ by $\Phi^{-1}H^iX=0$ for all $i$.
The distinction is between a single annihilator for the object and
annihilators depending on cohomological degree:
\[
\begin{array}{c|c|c}
 &\text{uniform annihilation}&\text{degreewise torsion}\\ \hline
 \mathbf D^-_{\mathrm{fg}}(R)&\mathcal G^-_{W_\Phi}=\mathcal S^-_\Phi
                         &\mathcal N^-_\Phi=\operatorname{PC}^-(\mathcal S^-_\Phi)\\
 \mathbf D^+_{\mathrm{fg}}(R)&\mathcal F^+_{W_\Phi}=\mathcal H_R\cap\mathcal S^+_\Phi
                         &\mathcal N^+_\Phi=\operatorname{PC}^+(\mathcal S^+_\Phi)\\
 \mathbf D^+(R\text{-}\mathrm{mod})
 &\mathcal F^{+,\mathrm{ft}}_{W_\Phi}=\mathcal H_R^{+,\mathrm{ft}}\cap\mathcal S^{+,\mathrm{ft}}_\Phi
 &\mathcal N^{+,\mathrm{ft}}_\Phi=\operatorname{PC}^{+,\mathrm{ft}}(\mathcal S^{+,\mathrm{ft}}_\Phi).
\end{array}
\]
The negative annihilator formula
\begin{equation}\label{eq:tensor-uniform-generation}
 \mathcal G^-_W=\{X:V(\operatorname{Ann}X)\subseteq W\}
\end{equation}
is \cite[Corollaries~2.13(3) and~2.14]{Matsui_2017}.
The positive formulas follow from
Proposition~\ref{thm:perfect-hom-classification}, since
$V(I)\subseteq W_\Phi$ exactly when $I\cap\Phi\ne\varnothing$.

For any family $\mathcal E$, the smallest Serre subcategory containing
$\mathcal A_\Phi$ and all $H^iE$ is
$\operatorname{Serre}(\mathcal A_\Phi\cup\{H^iE\})$.
Its support determines the Postnikov closure of
$\operatorname{Thick}^{\odot}(\mathcal S^\pm_\Phi\cup\mathcal E)$, namely
$\mathbf D^\pm_{\mathrm{fg}}(R)_{W_\Phi\cup\operatorname{Supp}\mathcal E}$.

\subsection{Relative comparisons and fixed support}
Let $F=\Phi^{-1}(-)$ and let $\mathcal E\subseteq\mathbf D^\pm_{\mathrm{fg}}(R)$.
Use $\operatorname{Thick}^{\diamond}$ for ordinary generation, or for
tensor generation in the negative category and Hom generation in the
positive category. Write
\[
\begin{aligned}
 \mathcal C&=\operatorname{Thick}^{\diamond}(\mathcal E),
 &\mathcal C_\Phi&=\operatorname{Thick}^{\diamond}(F\mathcal E),\\
 \mathcal J_{\mathrm{alg}}&=\operatorname{Thick}^{\diamond}(\mathcal C\cup\mathcal S^\pm_\Phi),
 &\mathcal J_{\mathrm{geo}}&=\operatorname{Thick}^{\diamond}(\mathcal C\cup\mathcal N^\pm_\Phi),\\
 \mathcal P&=\mathbf D^\pm_{\mathrm{fg}}(R)_{W_\Phi\cup\operatorname{Supp}\mathcal E}.
\end{aligned}
\]

\begin{proposition}\label{thm:relative-localization-reconstruction}
\label{thm:abelian-derived-saturation}
The relative comparison is
\begin{equation}\label{eq:relative-localization-ladder}
 \frac{\mathbf D^\pm_{\mathrm{fg}}(R)}{\mathcal J_{\mathrm{alg}}}
 \longrightarrow
 \frac{\mathbf D^\pm_{\mathrm{fg}}(\Phi^{-1}R)}{\mathcal C_\Phi}
 \longrightarrow
 \frac{\mathbf D^\pm_{\mathrm{fg}}(\Phi^{-1}R)}{\operatorname{PC}^\pm(F\mathcal E)}.
\end{equation}
The two kernels are respectively
$\mathcal J_{\mathrm{geo}}/\mathcal J_{\mathrm{alg}}$ and
$\operatorname{PC}^\pm(F\mathcal E)/\mathcal C_\Phi$.
They fit into the Verdier localization sequence
\begin{equation}\label{eq:relative-defect-sequence}
 \frac{\mathcal J_{\mathrm{geo}}}{\mathcal J_{\mathrm{alg}}}
 \longrightarrow\frac{\mathcal P}{\mathcal J_{\mathrm{alg}}}
 \longrightarrow\frac{\operatorname{PC}^\pm(F\mathcal E)}{\mathcal C_\Phi}.
\end{equation}
The same assertions hold internally in $\mathbf D^+(R\text{-}\mathrm{mod})$.
\end{proposition}
\begin{proof}
The heart consists of finite modules, so its wide subcategories are
Serre \cite[Corollary~3.2]{Takahashi}. Serre quotient correspondence and
Proposition~\ref{thm:postnikov-support} give
\[
 \operatorname{PC}^{\pm}(\mathcal J_{\mathrm{alg}})
  =\operatorname{PC}^{\pm}(\mathcal J_{\mathrm{geo}})=\mathcal P,
 \qquad F^{-1}\operatorname{PC}^{\pm}(F\mathcal E)=\mathcal P.
\]
Proposition~\ref{prop:central-comparison} identifies $F$ as the Verdier
quotient by $\mathcal N^\pm_\Phi$. Acting complexes lift by
\cite[Proposition~6.18]{sanders2026balmer}, and
\eqref{eq:finite-term-hom-diagonal} identifies the localized Hom actions.
Thus the middle term in \eqref{eq:relative-localization-ladder} is
$\mathbf D^\pm_{\mathrm{fg}}(R)/\mathcal J_{\mathrm{geo}}$.
Verdier transitivity \cite[Tag~05RJ]{stacks-project} gives both kernels
and \eqref{eq:relative-defect-sequence}. The same argument uses the
internal positive localization of
Proposition~\ref{prop:geometric-localization-plus} for finite terms.
\end{proof}

\begin{example}\label{ex:relative-two-defects}
Let $R=k[x,y]$, invert $x$, and take $\mathcal E=\{K(y)\}$.
Then \eqref{eq:relative-localization-ladder} is
\[
 \mathbf D^\pm_{\mathrm{fg}}(R)/\mathcal J_{\mathrm{alg}}
 \longrightarrow
 \mathbf D^\pm_{\mathrm{fg}}(R_x)/\operatorname{Thick}^{\diamond}(K(y)_x)
 \longrightarrow\mathbf D^\pm_{\mathrm{fg}}(R_{xy}).
\]
Both kernels are nonzero. Set
$U_f^-=\bigoplus_{n\ge1}K(f^n)[n]$ and $U_f^+=U_f^-\odot R$.
The image of $U_x^\pm$ belongs to the first kernel, and the image of
$(U_y^\pm)_x$ to the second.
Every object in $\mathcal J_{\mathrm{alg}}$ is killed by some $x^ay^b$.
If $U_x^\pm$ belonged to $\mathcal J_{\mathrm{alg}}$, localization
at $(x)$ and the perfect direct summands $K(x^n)$, or their shifted
duals in the positive case, would force $a\ge n$ for every $n$,
a contradiction. After inverting $x$, the same argument with
$y$ excludes $(U_y^\pm)_x$ from
$\operatorname{Thick}^{\diamond}(K(y)_x)$.
For ordinary thick generation this subcategory is
$\mathbf D^b_{\mathrm{fg}}(R_x)_{V(y)}$, so the second kernel is
$\mathbf D^\pm_{\mathrm{fg}}(R_x)_{V(y)}/\mathbf D^b_{\mathrm{fg}}(R_x)_{V(y)}$.
\end{example}

\begin{corollary}\label{cor:full-growth-common-closure}
\label{thm:growth-in-completion-kernel}
Fix a ring application $u:R\to V$ of
Theorem~\ref{thm:valuation-growth-retraction}, with
$u(r)=\varepsilon\pi^c$, $c\ge1$.
For every nonzero radical ideal $\mathcal I\subseteq\mathbf D^-_{\mathrm{fg}}(V)$
consisting of torsion objects, $\mathcal J=\Sigma_r(\mathcal I)$ satisfies
\begin{gather}
 \operatorname{Supp}\mathcal J=V(r),\qquad
 \operatorname{PC}^-(\mathcal J)=\mathcal N^-_r,
       \label{eq:growth-common-completion}\\
 \mathcal J\cap\mathbf D^b_{\mathrm{fg}}(R)=\mathbf D^b_{\mathrm{fg}}(R)_{V(r)},
 \qquad\mathcal J\cap\mathbf D^{\mathrm{perf}}(R)=\mathbf D^{\mathrm{perf}}(R)_{V(r)}.
       \label{eq:growth-common-perfect-part}
\end{gather}
The assignment $\mathcal I\mapsto\mathcal J$ preserves and reflects
all inclusions. In particular $J_{\lceil(n+1)^\alpha\rceil}\subseteq
J_{\lceil(n+1)^\beta\rceil}$ exactly when $\alpha\le\beta$.
If $R$ is not artinian, $r$ can be chosen so that, for
$f_j(n)=\lceil(n+1)^{1/j}\rceil$ and $\mathcal J_j=J_{f_j}$,
\begin{equation}\label{eq:descending-radical-growth}
 \sqrt{\mathcal S^-_r}\subsetneq\mathcal J_{j+1}
       \subsetneq\mathcal J_j\subsetneq\mathcal N^-_r\qquad(j\ge1).
\end{equation}
Their images form a strict chain inside $\ker\overline L^-_r$;
all bounded standard truncations of their generators vanish there.
\end{corollary}
\begin{proof}
Every positive $f$ has $\operatorname{Supp}K_f(r)=V(r)$.
The DVR classification expresses $\mathcal I$ as a nonempty join of
the principal radical ideals $\beta(K_f(\pi))$, so $\Sigma_r(\mathcal I)$ has this support.
Its intersections with the bounded and perfect subcategories are
given by \cite[Proposition~2.9 and Corollary~2.11]{Matsui_2017};
Proposition~\ref{thm:postnikov-support} computes its Postnikov closure.
Order reflection is $F_*\Sigma_r=\operatorname{id}$.
The comparison $J_f\subseteq J_g\Longleftrightarrow f\le_\mu g$
distinguishes the powers $\alpha$ and gives the strict chain.
For $b(n)=1$,
\[
 K_b(r)\simeq\operatorname{Cone}(r)\otimes_R^{\mathbf L}
                  \bigoplus_{n\ge0}R[n],\qquad
 \operatorname{Cone}(r)\mid K_b(r),
\]
so $J_b=\sqrt{\mathcal S^-_r}$.
The function $2^n$ gives an object of $\mathcal N^-_r$ outside every
$\mathcal J_j$.
For a nonartinian ring choose $u,r$ as in
Corollary~\ref{thm:noetherian-balmer-spectrum}.
The ideals contain $\mathcal S^-_r$, so their strictness survives the
quotient. Their bounded truncations are in $\mathcal G^-_{V(r)}=\mathcal S^-_r$.
\end{proof}

\begin{proposition}\label{cor:hom-growth-completion-defect}
For $u,r$ as above, consider
\[
 E^{+,\mathrm{ft}}=\bigoplus_{n\ge0}\operatorname{Cone}(r^{n+1})^\vee[-n]
       \in\mathbf D^+(R\text{-}\mathrm{mod}),\qquad
 E^+=j_R(E^{+,\mathrm{ft}}).
\]
Then $E^+\in\mathcal H_R\cap\mathcal N^+_r\setminus\mathcal S^+_r$
and $E^{+,\mathrm{ft}}\in\mathcal H_R^{+,\mathrm{ft}}\cap
\mathcal N^{+,\mathrm{ft}}_r\setminus\mathcal S^{+,\mathrm{ft}}_r$.
Both have support $V(r)$, and
\[
 \mathcal F^+_{V(r)}\subsetneq\mathcal H_R\cap\mathcal N^+_r,
 \qquad
 \mathcal F^{+,\mathrm{ft}}_{V(r)}
     \subsetneq\mathcal H_R^{+,\mathrm{ft}}\cap\mathcal N^{+,\mathrm{ft}}_r.
\]
The two subcategories $\mathcal F^+_{V(r)}$ and
$\mathcal F^{+,\mathrm{ft}}_{V(r)}$ have Postnikov closure equal to the
corresponding $\mathcal N_r$.
\end{proposition}
\begin{proof}
The sum has at most two finite free terms in each degree.
It is the internal Hom $K_f(r)\odot R$, with $f(n)=n+1$, so it belongs to
$\mathcal H_R^{+,\mathrm{ft}}$, and its image belongs to $\mathcal H_R$.
It becomes acyclic after inverting $r$ and has support $V(r)$.
Every summand is an internal perfect retract. Its base change has
\[
 H^{n+1}\!\left(V\otimes_R^{\mathbf L}
                   \operatorname{Cone}(r^{n+1})^\vee[-n]\right)
       =V/\pi^{c(n+1)}V.
\]
Thus $r^b\operatorname{id}_{E^+}\ne0$ for $n+1>b$, and the same
holds internally before applying $j_R$. This uses base change only
on the perfect retracts. The uniform-annihilator descriptions exclude
both objects from their $\mathcal F_{V(r)}$; the closure statement follows
from Proposition~\ref{thm:postnikov-support}.
\end{proof}

Consequently, $R$ is artinian exactly when every positive thick
Hom-submodule is Postnikov-complete, in either positive realization,
and exactly when, for every $\Phi$, the three comparisons
$\overline L^-_\Phi$, $\overline L^+_\Phi$ and
$\overline L^{+,\mathrm{ft}}_\Phi$ are equivalences.
The complete fourteen-condition statement, with its relative version
and proofs, is given in Appendix~\ref{subsec:artinian-equivalences}.

\phantomsection
\label{obs:group-postnikov-collapse}
\label{cor:group-postnikov-collapse}
For a finite group $G$ and field $k$, let $\mathcal T$ be any of
\[
 \mathbf D^b(kG\text{-}\mathrm{mod}),\quad
 \mathbf D^-(kG\text{-}\mathrm{mod}),\quad
 \mathbf D^+(kG\text{-}\mathrm{mod}),\quad
 \mathbf D^+_{\mathrm{fg}}(kG\text{-}\mathrm{Mod}).
\]
Let $\mathcal K=\mathbf D^b(kG\text{-}\mathrm{mod})$ act on $\mathcal T$
by $A\odot X=A\otimes_kX$. For every nonzero thick $\mathcal K$-submodule $\mathcal I$,
\[
 \operatorname{PC}(\mathcal I)=\mathcal T,\qquad
 \operatorname{Thick}^{\odot}(kG)=\mathbf K^b(kG\text{-}\mathrm{proj}).
\]
Indeed, if $0\ne H^nX$, the projection formula gives
\[
 H^n(kG\odot X)\cong(kG)^{\dim_kH^nX}
\]
\cite[Corollary~4.3.8(4)]{webb2016representations}.
The smallest wide subcategory containing all cohomology modules of
objects of $\mathcal I$ therefore contains $kG$. Finite presentations
then give every finite $kG$-module. Apply
Proposition~\ref{thm:abstract-postnikov-heart}.
If $\operatorname{char}k\mid|G|$, the perfect subcategory is proper
and the bounded category is rigid, so Postnikov closure is not
injective on its radical tensor ideals. This is the bounded modular
support setting of \cite[Theorem~11.4]{benson2011modular}.
In contrast, the full half-bounded category has nonnoetherian spectrum
by Theorem~\ref{thm:all-group-noetherian-ring}.

\subsection{A question on half-bounded spectra}

Let $\mathcal U$ be a symmetric tensor triangulated category with
coproducts preserved by tensor product and a nondegenerate
$t$-structure whose cohomology commutes with coproducts. Suppose
that its heart is a locally noetherian Grothendieck category,
$\mathbf1\in\mathcal U^{\le0}$, and
$\mathcal U^{\le0}\otimes\mathcal U^{\le0}\subseteq\mathcal U^{\le0}$.
Consider the full subcategory
\[
 \mathcal T=\{E\in\mathcal U\mid H^iE=0\ (i\gg0),\quad
                H^iE\text{ is noetherian for every }i\}.
\]
Assume that $\mathcal T$ is essentially small, contains $\mathbf1$,
and is closed under tensor product, and that $\mathcal U$ is
generated under triangles and coproducts by dualizable objects
of $\mathcal T$.

\begin{conjecture}\label{conj:postnikov-noetherianity}
Under these hypotheses,
\[
 \operatorname{Spc}(\mathcal T)\text{ is noetherian}
 \quad\Longleftrightarrow\quad
 \mathcal I=\operatorname{PC}^-(\mathcal I)
 \text{ for every radical thick tensor ideal }\mathcal I.
\]
\end{conjecture}

By Proposition~\ref{thm:abstract-postnikov-heart}, the condition on the
right means that each radical ideal is determined by its cohomology.
The full half-bounded category is essential: the bounded modular
category has noetherian spectrum, while its nonzero proper radical
ideals have whole-category Postnikov closure
\cite[Theorem~11.4]{benson2011modular}.

The conjecture holds in the ring, scheme, arithmetic and group
families considered above. In the noetherian cases the claims follow
from the artinian support classification, averaging, and rational
splitting. For nonartinian coefficients or positive-dimensional
schemes, the constructions above give radical thick tensor ideals
that are not Postnikov-complete. Their inverse images under the
$t$-exact forgetful or affine restriction functors give the
corresponding examples. For nonrational $R=S^{-1}\mathbb Z$, choose a noninvertible prime $p$. Then
\begin{gather*}
 E=\bigvee_{n\ge0}\operatorname{Cone}(p^{n+1}:\mathbb S[S^{-1}]\to\mathbb S[S^{-1}])[n]
      \in\operatorname{PC}^-(\mathcal S_p)\setminus\sqrt{\mathcal S_p},\\
 \pi_{rn}(\operatorname{Cone}(p^{n+1})^{\wedge r}[rn])=R/p^{n+1}R
 \qquad(r\ge1).
\end{gather*}
For modular groups, Proposition~\ref{prop:local-modular-growth} gives
proper radical ideals whose Postnikov closure contains $RG$ and
hence every finite module, by the preceding observation.

\clearpage
\appendix
\section{Central localization}
\label{app:central-localization}\label{sec:central-localization}

We record the module forms and image formulas for the central localization
used in Proposition~\ref{prop:minimal-central-localization}.

\subsection{Abelian categories}\label{sec:algebraic-abelian}

Let $\mathcal A$ be symmetric monoidal abelian and let $\mathcal B$
be an abelian $\mathcal A$-module, with additive right exact tensor
product and action $\odot$. We write
\[
 C=\operatorname{End}_{\mathcal A}(\mathbf1),\qquad
 \mathcal S=\{X\in\mathcal B:\phi_X=0\text{ for some }\phi\in\Phi\},
 \qquad\phi_X=\phi\odot X,
\]
where $\Phi\subseteq C$ is multiplicative. Write
$\mathcal U*\mathcal V$ for extensions of objects of $\mathcal V$
by objects of $\mathcal U$.

\begin{lemma}\label{lem:scalar-torsion-serre}
$\mathcal S$ is a Serre $\mathcal A$-submodule.
\end{lemma}
\begin{proof}
Naturality handles subobjects, quotients and the action. In
$0\to X\xrightarrow{i}Y\xrightarrow{p}Z\to0$,
\[
 s_X=t_Z=0\ \Longrightarrow\ p t_Y=0
 \ \Longrightarrow\ t_Y=ih
 \ \Longrightarrow\ (st)_Y=i s_Xh=0.
\]
\end{proof}

\begin{proposition}
\label{prop:abelian-localization-hom}
\label{prop:localization-serre-generation}
\label{prop:preimage-serre-localization}
\label{thm:localization-commutes-serre-quotients}
For $\pi:\mathcal B\to\mathcal B/\mathcal S$,
\[
 \Phi^{-1}\operatorname{Hom}_{\mathcal B}(X,Y)
 \xrightarrow{\sim}\operatorname{Hom}_{\mathcal B/\mathcal S}(\pi X,\pi Y),
 \qquad f/s\longmapsto\pi(s_Y)^{-1}\pi f.
\]
The action descends. If $\mathcal C=\operatorname{Serre}^{\odot}(\mathcal E)$
and $q:\mathcal B\to\mathcal B/\mathcal C$, then
\begin{gather*}
 \pi\mathcal C=\operatorname{Serre}^{\odot}(\pi\mathcal E),\qquad
 \pi^{-1}(\pi\mathcal C)=\mathcal S*\mathcal C
 =\mathcal C*\mathcal S
 =\operatorname{Serre}^{\odot}(\mathcal S\cup\mathcal E),\\
 \mathcal S_{\mathcal B/\mathcal C}(\Phi)=q\mathcal S,\qquad
 \frac{\mathcal B/\mathcal S}
      {\operatorname{Serre}^{\odot}(\pi\mathcal E)}
 \simeq
 \frac{\mathcal B/\mathcal C}{q\mathcal S}.
\end{gather*}
\end{proposition}
\begin{proof}
The kernel and cokernel of $s_Y$ are killed by $s$. Conversely,
if $u:X\to Y$ has kernel and cokernel in $\mathcal S$, its image
factorization gives $v:Y\to X$ and $s\in\Phi$ with
$vu=s_X$, $uv=s_Y$. Thus Gabriel fractions are scalar fractions,
and
\[
 \pi f=0\iff\operatorname{im}f\in\mathcal S
           \iff tf=0\text{ for some }t\in\Phi.
\]
Each functor defined by the action preserves these scalar inverses,
so the action descends.

In an extension of $S\in\mathcal S$ by $C\in\mathcal C$, or of
$C$ by $S$, the image of a scalar killing $S$ belongs to
$\mathcal C$, while its kernel and cokernel belong to $\mathcal S$.
Hence $\mathcal S*\mathcal C=\mathcal C*\mathcal S$. Subobjects
and quotients preserve this class; exchanging its two middle
factors proves extension closure. Right exactness gives action
stability. Scalar fractions now identify this class with
$\pi^{-1}(\pi\mathcal C)$, and Serre correspondence gives the
first two formulas.

Finally, $q(s_X)=0$ implies $\operatorname{im}s_X\in\mathcal C$,
so $qX\simeq q\ker s_X$ with $\ker s_X\in\mathcal S$.
Thus the torsion image is $q\mathcal S$. Both iterated quotients
are the Gabriel quotient by $\operatorname{Serre}^{\odot}
(\mathcal C\cup\mathcal S)$.
\end{proof}

\begin{corollary}\label{cor:finite-module-central-localization}
Let $R$ be commutative noetherian and $\Phi\subseteq R$ multiplicative.
Then
\[
 (R\text{-}\mathrm{mod})/
 \{M:sM=0\text{ for some }s\in\Phi\}
 \simeq(\Phi^{-1}R)\text{-}\mathrm{mod}.
\]
\end{corollary}
\begin{proof}
Finite presentation gives the localized Hom formula; clearing the
finitely many entries in a presentation matrix gives density.
\end{proof}

\subsection{Triangulated categories}\label{sec:algebraic-triangulated}

Let $\mathcal K$ be an essentially small exact triangulated module
over an essentially small tensor triangulated category $\mathcal T$,
with action $\odot$ as in \cite[Definition~3.2(1)--(3)]{Stevenson}.
Let $C=\operatorname{End}_{\mathcal T}(\mathbf1)$ and let
$\Phi\subseteq C$ be multiplicative. Write
$\mathcal S=\{X:\phi_X=0\text{ for some }\phi\in\Phi\}$.
Here $\mathcal U*\mathcal V$ denotes extension by triangles.

\begin{lemma}\label{lem:triangulated-uniform-torsion}
$\mathcal S$ is a thick submodule and
$\phi^2\operatorname{id}_{\operatorname{Cone}(\phi_X)}=0$.
\end{lemma}
\begin{proof}
The product formula for annihilators and the cone calculation in the proof of
Proposition~\ref{prop:minimal-central-localization} apply objectwise
to the module $\mathcal K$.
\end{proof}

\begin{proposition}
\label{prop:triangulated-localization-hom}
\label{prop:triangulated-localization-generation}
\label{thm:localization-commutes-triangulated-quotients}
For $\pi:\mathcal K\to\mathcal K/\mathcal S$,
\[
 \Phi^{-1}\operatorname{Hom}_{\mathcal K}(X,Y)
 \xrightarrow{\sim}\operatorname{Hom}_{\mathcal K/\mathcal S}(\pi X,\pi Y),
 \qquad f/s\longmapsto\pi(s_Y)^{-1}\pi f.
\]
The action descends. For $\mathcal C=\operatorname{Thick}^{\odot}(\mathcal E)$
and $q:\mathcal K\to\mathcal K/\mathcal C$,
\begin{gather*}
 \operatorname{Thick}^{\odot}(\pi\mathcal E)
       =\operatorname{smd}(\pi\mathcal C),\\
 \pi^{-1}\operatorname{Thick}^{\odot}(\pi\mathcal E)
       =\operatorname{smd}(\mathcal S*\mathcal C)
       =\operatorname{smd}(\mathcal C*\mathcal S)
       =\operatorname{Thick}^{\odot}(\mathcal S\cup\mathcal E),\\
 \mathcal S_{\mathcal K/\mathcal C}(\Phi)
       =\operatorname{smd}(q\mathcal S),\qquad
 \frac{\mathcal K/\mathcal S}
      {\operatorname{Thick}^{\odot}(\pi\mathcal E)}
       \simeq
 \frac{\mathcal K/\mathcal C}{\operatorname{smd}(q\mathcal S)}.
\end{gather*}
The formulas also hold for ordinary thick generation, with action
compatibility omitted if $\mathcal C$ is not a submodule.
The ordinary formulas hold for any degree-zero central action of $C$,
by the same scalar-fraction argument.
\end{proposition}
\begin{proof}
The fraction formula is degree-zero central localization
\cite[Corollary~3.10]{balmer2010spectra}. For modules the same
roof argument applies: if $\operatorname{Cone}(u)$ is killed by
$s\in\Phi$, the triangle of $u$ gives $v$ with $uv=s$, and
the roof $(u,f)$ represents $fv/s$. The Verdier zero criterion
gives injectivity. The action descends because $\mathcal S$ is a
submodule.

Since $\operatorname{Cone}(f/s)\simeq\pi\operatorname{Cone}(f)$,
$\pi\mathcal C$ is triangulated. An isomorphism
$\pi X\simeq\pi C$ and its inverse are scalar fractions, giving
$\pi^{-1}(\pi\mathcal C)=\mathcal S*\mathcal C
=\mathcal C*\mathcal S$. Lifting a complement of any summand
of $\pi C$ proves the two generation formulas.

If $q(\phi_X)=0$, then
\[
 q\operatorname{Cone}(\phi_X)\simeq qX\oplus qX[1],
 \qquad \operatorname{Cone}(\phi_X)\in\mathcal S.
\]
This proves the torsion-image formula. Both successive quotients
are $\mathcal K/\operatorname{Thick}^{\odot}(\mathcal C\cup\mathcal S)$
by the Verdier universal property \cite[Tag~05RJ]{stacks-project}.
The ordinary case uses the same argument.
\end{proof}

\subsection{Actions on the degreewise kernel}

Use $\mathcal S_\Phi$, $\mathcal N_\Phi$ and $\mathcal T^\pm$
from Proposition~\ref{prop:central-comparison}.

\begin{proposition}
\label{prop:abstract-action-degreewise-kernel}
Suppose that $\mathcal T$ is an exact triangulated module over a tensor
triangulated category $\mathcal K$, that its scalars arise from the
unit, and that every object $A\in\mathcal K$ has an integer $d_A$ such
that, for every $n$,
\[
 \begin{cases}
 A\odot\mathcal T^{\leq n}\subseteq\mathcal T^{\leq n+d_A},
       &\mathcal T=\mathcal T^-,\\
 A\odot\mathcal T^{\geq n}\subseteq\mathcal T^{\geq n+d_A},
       &\mathcal T=\mathcal T^+.
 \end{cases}
\]
Then both $\mathcal S_\Phi$ and $\mathcal N_\Phi$ are thick
$\mathcal K$-submodules, so all comparisons in
Proposition~\ref{prop:central-comparison} carry the induced actions.
\end{proposition}
\begin{proof}
Scalar annihilators are preserved by the action.  For
$X\in\mathcal N_\Phi$, set $m=i-d_A$.  In the negative case the
truncation triangle identifies
\[
 H^i(A\odot X)\cong H^i(A\odot\tau^{\geq m}X),
\]
since $A\odot\tau^{\leq m-1}X\in\mathcal T^{\leq i-1}$.
The truncation is bounded and has a uniform annihilator, so this
cohomology lies in $\mathcal A_\Phi$.  The positive case uses
$H^i(A\odot X)\cong H^i(A\odot\tau^{\leq m}X)$, since the remaining
term lies in $\mathcal T^{\geq i+1}$.  Thus $\mathcal N_\Phi$ is a
submodule.
\end{proof}

\begin{remark}\label{rem:central-nondegenerate-factorization}
In Proposition~\ref{prop:central-comparison}, the quotient
$\mathcal Q_\Phi^\pm/\mathcal D_\Phi^\pm
\simeq\mathcal T^\pm/\mathcal N_\Phi^\pm$ has a nondegenerate
$t$-structure with the same heart. For the minus sign, a roof
\[
 X\xleftarrow{v} Z\xrightarrow{f}Y,\qquad
 X\in(\mathcal Q_\Phi^-)^{\le0},\quad
 Y\in(\mathcal Q_\Phi^-)^{\ge1},\quad
 \operatorname{Cone}(v)\in\mathcal D_\Phi^-
\]
has $Z\in(\mathcal Q_\Phi^-)^{\le0}$ and hence $f=0$.
Truncation triangles therefore descend. Truncating roofs identifies
the heart; the zero-cohomology objects upstairs are precisely
$\mathcal D_\Phi^-$, proving nondegeneracy. The plus sign is dual.

Consequently, a $t$-exact $F:\mathcal Q_\Phi^\pm\to\mathcal V$,
with $\bigcap_n\mathcal V^{\le n}=0$ for the minus sign or
$\bigcap_n\mathcal V^{\ge n}=0$ for the plus sign, satisfies
\[
 F(\mathcal D_\Phi^\pm)=0,\qquad
 F=\overline F\lambda,\qquad
 \mathcal Q_\Phi^\pm\xrightarrow{\lambda}
 \mathcal T^\pm/\mathcal N_\Phi^\pm
       \xrightarrow{\overline F}\mathcal V.
\]
The $t$-exact factor $\overline F$ is unique up to natural isomorphism
by the Verdier universal property \cite[Tag~05RJ]{stacks-project}.
\end{remark}

\begin{example}\label{ex:triangulated-smd-essential}
Let $R=k[x,y]/(xy)$, $s=x+y$, and $\Phi=\{s^n:n\ge0\}$.
In $\mathbf D^b_{\mathrm{fg}}(R)$, let
$P=[R\xrightarrow{y}R]$ in degrees $-1,0$ and
$\mathcal C=\operatorname{Thick}^{\odot}(P)$ for the perfect
tensor action. Then
\[
 \pi(R/(y))\in\operatorname{smd}(\pi\mathcal C)\setminus\pi\mathcal C.
\]
The injection $i:R/(y)\to R$ given by multiplication by $x$ and
the quotient $p:R\to R/(y)$ satisfy $ip=x_R$ and
$pi=s_{R/(y)}$. Thus $\pi i/s$ splits $\pi p$, and the analogous
maps for $R/(x)$ give
\[
 \pi R\simeq\pi(R/(y))\oplus\pi(R/(x)),\qquad
 \pi P\simeq\pi(R/(y))\oplus\pi(R/(y))[1].
\]
Every object of $\mathcal C$ is perfect and has constant Euler
rank on the connected space $\operatorname{Spec}R$. But
\[
 R_s\simeq k[x,x^{-1}]\times k[y,y^{-1}],\qquad
 (R/(y))_s\simeq(k[x,x^{-1}],0)
\]
has ranks $1$ and $0$. Hence it cannot be the geometric
localization of an object of $\mathcal C$, proving the claim.
\end{example}
\section{Positive localization and Hom generation}
\label{app:positive-details}\label{sec:geometric-plus}

\subsection{A geometric lifting criterion}
\begin{proposition}
\label{prop:abstract-positive-geometric-lifting}
Let $L:\mathcal U\rightleftarrows\mathcal V:j$ be a $C$-linear
adjunction of exact functors, with $j$ fully faithful, $L$ $t$-exact, and $j$
preserving bounded-below objects.  Assume that both $t$-structures are
nondegenerate, both categories admit countable coproducts, and their
cohomology commutes with those coproducts.  Assume also that their
hearts have exact sequential colimits.  Let
$\mathcal A\subseteq\mathcal U^\heartsuit$ and
$\mathcal B\subseteq\mathcal V^\heartsuit$ be Serre subcategories,
and write
\[
 \mathcal U^+_{\mathcal A}=\{X\in\mathcal U^+\mid
                                       H^iX\in\mathcal A\ (\forall i)\},
 \qquad
 \mathcal V^+_{\mathcal B}=\{Y\in\mathcal V^+\mid
                                       H^iY\in\mathcal B\ (\forall i)\}.
\]
Suppose that $\Phi$ acts invertibly on $\mathcal V$, that $L$ carries
$\mathcal A$ into $\mathcal B$, and that
\begin{enumerate}
\item every $B\in\mathcal B$ is isomorphic to $LM$ for some
$M\in\mathcal A$;
\item for $M\in\mathcal A$, $E\in\mathcal U^+$, and $r\in\mathbb Z$,
the canonical map
\[
 \Phi^{-1}\operatorname{Hom}_{\mathcal U}(M[-r],E)
       \longrightarrow\operatorname{Hom}_{\mathcal V}(LM[-r],LE)
\]
is surjective.
\end{enumerate}
Then $L:\mathcal U^+_{\mathcal A}\to\mathcal V^+_{\mathcal B}$ is a
Verdier localization.  Its kernel is determined by the kernel of
$L:\mathcal A\to\mathcal B$ degree by degree.
\end{proposition}

\begin{proof}
We first prove an approximation assertion.  Given $E\in\mathcal U^+$
with $LE\in\mathcal V^+_{\mathcal B}$, choose $m$ with
$E\in\mathcal U^{\geq m}$.  Starting with $Z_{m-1}=0$, construct
compatible $Z_n\to E$ such that
\[
 Z_n\in\mathcal U^{[m,n]}_{\mathcal A},\qquad
 L\operatorname{Cone}(Z_n\to E)\in\mathcal V^{\geq n+1}.
\]
For $C_{n-1}=\operatorname{Cone}(Z_{n-1}\to E)$, choose
$M_n\in\mathcal A$ lifting $H^n(LC_{n-1})$.
The lowest truncation map $LM_n[-n]\to LC_{n-1}$ lifts after
multiplication by a scalar in $\Phi$, by (2).
The homotopy pullback then extends $Z_{n-1}\to E$ to $Z_n\to E$
and gives a triangle
\[
 Z_{n-1}\longrightarrow Z_n\longrightarrow M_n[-n]
                                      \longrightarrow Z_{n-1}[1].
\]
The maps $H^iZ_{n-1}\to H^iZ_n$ are isomorphisms for $i<n$.
The telescope $Z=\operatorname{hocolim}_nZ_n$ therefore has
cohomology in $\mathcal A$ and is bounded below, by nondegeneracy.
Compatibility gives $Z\to E$; exactness and preservation of
coproducts by $L$ identify its localized cohomology with that of
$E$ in each degree.  Nondegeneracy on the target makes
$LZ\to LE$ an isomorphism.

Apply this assertion to $jY$ to obtain density.  For
$f:LX\to LY$, form in $\mathcal U^+$ the homotopy pullback of
$(j(f)\eta_X,-\eta_Y):X\oplus Y\to jLY$.
Its projection to $X$ becomes invertible under $L$ and the other
projection represents $f$.  Approximate this middle object by
$Z\in\mathcal U^+_{\mathcal A}$ as above.  This gives a roof in the
restricted category and proves fullness of the functor from its
Verdier quotient by $\ker L$.  That functor has zero kernel, hence
is faithful by the standard full exact functor criterion.
\end{proof}

\subsection{Localization of the two bounded-below derived categories}
\begin{proposition}\label{prop:geometric-localization-plus}
\label{prop:geometric-localization-finite-terms}\label{loc dense}
For a commutative noetherian ring $R$ and multiplicative $\Phi$,
both functors
\[
 \mathbf D^+_{\mathrm{fg}}(R)\longrightarrow
     \mathbf D^+_{\mathrm{fg}}(\Phi^{-1}R),\qquad
 \mathbf D^+(R\text{-}\mathrm{mod})\longrightarrow
     \mathbf D^+(\Phi^{-1}R\text{-}\mathrm{mod})
\]
are Verdier localizations. Their kernels consist of the objects whose
cohomology becomes zero degree by degree after localization.
\end{proposition}
\begin{proof}
For the first functor use flat localization and restriction of scalars
on $\mathbf D(R)$. The ambient $t$-structures are nondegenerate, their
cohomology commutes with coproducts, and their hearts have exact
sequential colimits. Finite presentation lifts each finite module over
$\Phi^{-1}R$ to a finite $R$-module. For a finite $M$, a bounded-below
$E$, and any $r$, the map
\[
 \Phi^{-1}\operatorname{Hom}_{\mathbf D(R)}(M[-r],E)
  \longrightarrow
 \operatorname{Hom}_{\mathbf D(\Phi^{-1}R)}
                     (\Phi^{-1}M[-r],\Phi^{-1}E)
\]
is surjective by \cite[Tag~0A6A(3)]{stacks-project}.
Proposition~\ref{prop:abstract-positive-geometric-lifting} applies and
gives the first Verdier localization, including object and morphism lifts.

For $\mathbf D^+(R\text{-}\mathrm{mod})$, use Miyachi's localization theorem
\cite[Theorem~3.2]{miyachi1991localization}; see also
\cite[Theorem~1.1]{wuzhou2022shortexact}.
With $\mathcal A_\Phi=\{M\in R\text{-}\mathrm{mod}:\Phi^{-1}M=0\}$,
finite presentation gives
$(R\text{-}\mathrm{mod})/\mathcal A_\Phi\simeq
(\Phi^{-1}R)\text{-}\mathrm{mod}$, and Miyachi gives
\[
 \mathbf D^+(R\text{-}\mathrm{mod})/
       \mathbf D^+_{\mathcal A_\Phi}(R\text{-}\mathrm{mod})
  \simeq\mathbf D^+(\Phi^{-1}R\text{-}\mathrm{mod}).
\]
The subcategory being quotiented out is exactly the stated kernel.
This proves the localization assertion for each of the two categories.
\end{proof}

\subsection{The artinian comparison}
\begin{remark}\label{rem:finite-term-comparison}
For a commutative noetherian ring $R$, the canonical comparison is an
equivalence
\[
 \mathbf D^-(R\text{-}\mathrm{mod})
 \xrightarrow{\sim}\mathbf D^-_{\mathrm{fg}}(R);
\]
see \cite[Proposition~2.1]{positselski2021unbounded}.
For artinian $R$, the canonical functor
\[
 \mathbf D^+(R\text{-}\mathrm{mod})
 \longrightarrow\mathbf D^+_{\mathrm{fg}}(R)
\]
is also an equivalence. Indeed, a finite module has a finite injective
hull: reduce to the artinian local factors and use
\cite[Tag~08YX]{stacks-project}. For a monomorphism $N\hookrightarrow M$
with $N$ finite, injectivity extends $N\hookrightarrow E(N)$ to
$M\to E(N)$. The comparison criterion of
\cite[Section~3, remark~(a) following Theorem~3.1]{positselski2021unbounded}
therefore applies. In the general noetherian case, the two positive
realizations remain separate throughout the paper.
\end{remark}

\subsection{Internal derived Hom}\label{app:internal-hom}
\begin{proof}[Details for Lemma~\ref{prop:derived-hom-module}]
Represent $A^\bullet$ by a bounded-above finite projective complex
$P^\bullet$, with $P^i=0$ for $i>b$, and represent $X^\bullet$
by a complex with $X^i=0$ for $i<a$. Then
\begin{equation}\label{eq:finite-term-hom-diagonal}
 \operatorname{Hom}_R^\bullet(P^\bullet,X^\bullet)^n
   =\bigoplus_{a-n\le i\le b}
            \operatorname{Hom}_R(P^i,X^{i+n}),
 \qquad \operatorname{Hom}_R^\bullet(P^\bullet,X^\bullet)^n=0
                  \quad(n<a-b).
\end{equation}
If every $X^i$ is finite, this is a bounded-below complex of finite
modules, so defines an object of $\mathbf D^+(R\text{-}\mathrm{mod})$. The complex
$P^\bullet$ is $K$-projective: Hom from it preserves
quasi-isomorphisms in $X^\bullet$, and changing its projective
resolution gives a homotopy equivalence. Thus this construction
descends internally to $\mathbf D^+(R\text{-}\mathrm{mod})$, without a full-faithfulness
assumption on $j_R$.

For the finite-cohomology category, a sufficiently low brutal cut
$P^{\ge m}$ leaves the Hom complex unchanged in degrees $n-1,n,n+1$.
Consequently its $n$th cohomology equals that of
$\mathbf R\operatorname{Hom}_R(P^{\ge m},X^\bullet)$, which is
finite because $P^{\ge m}$ is perfect. This proves well-definedness
in both categories.

Exactness follows from the Hom complexes and mapping cones.
For bounded-above finite projective $P^\bullet,Q^\bullet$, currying gives
\[
 \operatorname{Hom}_R^\bullet(P^\bullet\otimes_R Q^\bullet,X^\bullet)
 \cong
 \operatorname{Hom}_R^\bullet
 (P^\bullet,\operatorname{Hom}_R^\bullet(Q^\bullet,X^\bullet)).
\]
All diagonals are finite by the same bounds. Iterated currying
proves the pentagon, and evaluation at $1\in R$ proves the unit
diagrams. Under the standard Hom differential, the shift comparison
in the first variable has sign $(-1)^{rn}$ on Hom degree $n$;
shifting the second variable by $s$ changes it by $(-1)^{rs}$.
This is the required suspension compatibility. Finally the same
complexes compute the images under $j_R$; quasi-isomorphisms and
standard truncations are detected by cohomology. Hence $j_R$ has
the asserted properties.
\end{proof}
\section{Extended Artinian criteria}
\label{app:extended-criteria}
For a fixed multiplicative $\Phi\subseteq R$, the conditions of
Corollary~\ref{thm:localized-noetherian-spectrum} are also equivalent
to the heart $(\Phi^{-1}R)\text{-}\mathrm{mod}$ being a length category
and the induced $t$-structure on
$\mathbf D^-_{\mathrm{fg}}(R)/\mathcal S^-_\Phi$ being nondegenerate.
This follows from Proposition~\ref{prop:central-comparison}.

\subsection{The fourteen Artinian criteria}\label{subsec:artinian-equivalences}
The equivalence of (1), (3), (5) and (7) below is due to
Matsui--Takahashi \cite[Theorem~6.5]{Matsui_2017}, using
\eqref{eq:postnikov-tame}.
\begin{corollary}\label{thm:artinian-equivalences}
For a commutative noetherian ring $R$, the following are equivalent:
\begin{enumerate}
 \item $R$ is artinian;
 \item $\operatorname{Spc}(\mathbf D^-_{\mathrm{fg}}(R))$ is noetherian;
 \item every thick tensor ideal of $\mathbf D^-_{\mathrm{fg}}(R)$
       is Postnikov-complete;
 \item every thick Hom-submodule of $\mathbf D^+_{\mathrm{fg}}(R)$
       is Postnikov-complete;
 \item $\operatorname{PC}^-$ is injective on radical thick tensor ideals
       of $\mathbf D^-_{\mathrm{fg}}(R)$;
 \item $\operatorname{PC}^+$ is injective on thick Hom-submodules
       of $\mathbf D^+_{\mathrm{fg}}(R)$;
 \item for every specialization-closed $W\subseteq\operatorname{Spec}R$,
       $\mathcal G^-_W=\mathbf D^-_{\mathrm{fg}}(R)_W$;
 \item for every specialization-closed $W\subseteq\operatorname{Spec}R$,
       $\mathcal F^+_W=\mathbf D^+_{\mathrm{fg}}(R)_W$;
 \item for every multiplicative subset $\Phi\subseteq R$,
       $\overline L^-_\Phi$ is an equivalence;
 \item for every multiplicative subset $\Phi\subseteq R$,
       $\overline L^+_\Phi$ is an equivalence.

 \item every thick Hom-submodule of $\mathbf D^+(R\text{-}\mathrm{mod})$
       is Postnikov-complete;
 \item $\operatorname{PC}^{+,\mathrm{ft}}$ is injective on thick
       Hom-submodules of $\mathbf D^+(R\text{-}\mathrm{mod})$;
 \item for every specialization-closed $W\subseteq\operatorname{Spec}R$,
       $\mathcal F^{+,\mathrm{ft}}_W=(\mathbf D^+(R\text{-}\mathrm{mod}))_W$;
 \item for every multiplicative subset $\Phi\subseteq R$,
       $\overline L^{+,\mathrm{ft}}_\Phi$ is an equivalence.
\end{enumerate}

When these conditions hold, ordinary support classifies all the
tensor ideals and Hom-submodules in (3), (4) and (11); every one of them is
generated by perfect complexes under its respective action.
\end{corollary}

\begin{proof}
By \eqref{eq:postnikov-tame} and
\cite[Theorem~6.5(2), (9), (11)]{Matsui_2017}, (1), (3) and (5) are
equivalent. The bounds
$\mathcal G^-_W\subseteq\mathcal J\subseteq\mathcal J^{\mathrm{tame}}$
for $W=\operatorname{Supp}\mathcal J$
\cite[Lemma~5.8]{Matsui_2017} give (7).
Corollary~\ref{thm:noetherian-balmer-spectrum} gives $(1)\Longleftrightarrow(2)$.

If $R$ is artinian, Corollary~\ref{prop:artinian-hom-submodules} and
Proposition~\ref{thm:postnikov-support} give (4), (6), (8), (11), (12)
and (13), as well as the support and perfect-generation assertions.
For the idempotent $e$ of $U_\Phi$,
\[
 \Phi^{-1}R\cong eR,\qquad
 \phi_0(1-e)=0\text{ for some }\phi_0\in\Phi,\qquad
 \mathcal N_\Phi=\mathcal S_\Phi=\{X:X\simeq(1-e)X\}.
\]
Thus localization is the projection onto the $eR$ factor, proving
(9), (10) and (14) in their respective categories.

If $R$ is not artinian, Corollary~\ref{thm:growth-in-completion-kernel}
gives $\mathcal S^-_r\subsetneq\mathcal N^-_r$, so (9) fails by
Proposition~\ref{prop:central-comparison}.
Use the single pair $E^{+,\mathrm{ft}},E^+$ of
Proposition~\ref{cor:hom-growth-completion-defect}.
In each positive realization,
\[
 \mathcal F_{V(r)}\subseteq\mathcal S_r\subsetneq\mathcal N_r,
 \qquad E\in\mathcal N_r\setminus\mathcal S_r,\qquad
 \operatorname{Supp}E=V(r).
\]
The internal or finite-cohomology Hom-submodule generated by $E$ is
different from $\mathcal F_{V(r)}$, but both have support $V(r)$ and
Postnikov closure $\mathcal N_r$. This disproves (4), (6), (8) and,
respectively, (11), (12), (13). The nonzero kernel classes of $E$
disprove (10) and (14). Every argument here is read in its own positive
category; an equivalence between the two positive categories is not used.
\end{proof}

\subsection{Relative criteria after algebraic localization}

Let $R$ be commutative noetherian, let
$\mathcal C^-=\operatorname{Thick}(\mathcal E)$ be a tensor ideal in
$\mathbf D^-_{\mathrm{fg}}(R)$, and let $\Phi\subseteq R$ be multiplicative.
Define thick subcategories in the two bounded-below categories by
\[
\begin{aligned}
 \mathcal C^+&=\operatorname{Thick}
   \{A\odot Y:A\in\mathcal E,\ Y\in\mathbf D^+_{\mathrm{fg}}(R)\},\\
 \mathcal C^{+,\mathrm{ft}}&=\operatorname{Thick}
   \{A\odot Y:A\in\mathcal E,\ Y\in\mathbf D^+(R\text{-}\mathrm{mod})\}.
\end{aligned}
\]
Let $\mathcal D^\alpha=\operatorname{Thick}
(\mathcal C^\alpha\cup\mathcal S^\alpha_\Phi)$ for
$\alpha=-,+,(+,\mathrm{ft})$, and
\[
 W'=\operatorname{Supp}\mathcal E\cup\bigcup_{\phi\in\Phi}V(\phi),
 \qquad U=\operatorname{Spec}R\setminus W'.
\]
For each nonclosed $y\in U$, choose a homomorphism
$u_y:R\to V_y$ as in Theorem~\ref{thm:central-relative-scheme}.
Require, for each $A\in\mathcal E$ and
$Y\in\mathbf D^+_{\mathrm{fg}}(R)$, integers $a,b\ge1$ such that
\begin{equation}\label{eq:relative-positive-trait-condition}
 \pi_y^a\operatorname{id}_{V_y\otimes_R^{\mathbf L}A}=0,
 \qquad
 \pi_y^b\operatorname{id}_{V_y\otimes_R^{\mathbf L}(A\odot Y)}=0.
\end{equation}
The second equality is taken in $\mathbf D(V_y)$; its exponent
may depend on $A,Y,y$.

\begin{corollary}\label{cor:relative-artinian-criteria}
All fourteen conditions of Corollary~\ref{thm:artinian-equivalences}
remain equivalent with these replacements: (1) says that every
point of $U$ is closed in $\operatorname{Spec}R$; the three categories
are replaced respectively by
\[
 \frac{\mathbf D^-_{\mathrm{fg}}(R)}{\mathcal D^-},\qquad
 \frac{\mathbf D^+_{\mathrm{fg}}(R)}{\mathcal D^+},\qquad
 \frac{\mathbf D^+(R\text{-}\mathrm{mod})}{\mathcal D^{+,\mathrm{ft}}};
\]
and (7), (8), (13) use $Z\supseteq W'$ and the images of perfect
complexes supported on $Z$. Each is the algebraic $\Phi$-localization
of the quotient by its corresponding $\mathcal C^\alpha$.
For every further multiplicative set $\Psi$, (9), (10), (14) use
the relative comparisons of
Proposition~\ref{thm:relative-localization-reconstruction}.
When these conditions hold, $U$ is finite and clopen, and the
three quotients are the corresponding derived categories of the
artinian direct factor $R_U$.
\end{corollary}

\begin{proof}
The tensor--Hom identity and ordinary thick generation give
\[
 B\odot(A\odot Y)\simeq(A\otimes_R^{\mathbf L}B)\odot Y,
 \qquad
 \mathcal D^-\odot\mathbf D^+_{\mathrm{fg}}(R)\subseteq\mathcal D^+,
 \qquad
 \mathbf D^-_{\mathrm{fg}}(R)\odot\mathcal D^+\subseteq\mathcal D^+.
\]
Scalar annihilators handle $\mathcal S_\Phi$ in these inclusions.
The same proof works internally for finite terms, and the
canonical functor to $\mathbf D(R)$ commutes with $\odot$.
Thus both quotient actions descend, and
\eqref{eq:relative-positive-trait-condition} also applies to
finite-term generators.

Write $q^\alpha$ for the three quotients. Matsui--Takahashi's
bounded-part theorem \cite[Proposition~2.9]{Matsui_2017} gives
$\mathcal D^-\cap\mathbf D^b_{\mathrm{fg}}(R)
=\mathbf D^b_{\mathrm{fg}}(R)_{W'}$.
For a bounded $M$ supported on $W'$, choose generators $\mathbf a$
of $\operatorname{Ann}(\operatorname{id}_M)$. Then
\[
 K(\mathbf a)\in\mathcal D^-,\qquad
 M\mid K(\mathbf a)\odot M\in\mathcal D^+.
\]
Support containment for $\odot$ proves the reverse bounded-part
inclusion; the internal argument is identical.
For each subcategory $\mathcal D^\alpha$, one bounded truncation of
any $X\in\mathcal D^\alpha$ lies in
$\mathbf D^b_{\mathrm{fg}}(R)_{W'}$; its triangle puts the other
truncation in $\mathcal D^\alpha$. Thus $\mathcal D^\alpha$ is
truncation-closed. Write $\mathcal A_{W'}$ for its Serre heart.
To check quotient orthogonality, consider a roof
$X\xleftarrow{u}Z\to Y$, where $X\in\mathcal T^{\le0}$,
$Y\in\mathcal T^{\ge1}$ and $C=\operatorname{Cone}(u)\in\mathcal D^\alpha$.
The connecting morphism induces
\[
 \tau^{\ge0}C[-1]\longrightarrow\tau^{\ge1}Z,
 \qquad
 \operatorname{Cone}\bigl(\tau^{\ge0}C[-1]
             \longrightarrow\tau^{\ge1}Z\bigr)
 \simeq I,
\]
where $I=\operatorname{im}(H^0X\to H^0C)\in\mathcal A_{W'}$.
Hence $\tau^{\ge1}Z\in\mathcal D^\alpha$, and the roof can be
replaced by $X\leftarrow\tau^{\le0}Z\to Y$, whose numerator is
zero. Truncation triangles now give the quotient $t$-structure.
For heart objects the same reductions replace a roof by a weak
isomorphism in $R\text{-}\mathrm{mod}$ with kernel and cokernel in
$\mathcal A_{W'}$; its heart is therefore
$(R\text{-}\mathrm{mod})/\mathcal A_{W'}$.
Proposition~\ref{thm:abstract-postnikov-heart}, Takahashi's classification,
and Serre quotient correspondence yield, for thick subcategories
$\mathcal J^\pm\supseteq\mathcal D^\pm$,
\begin{equation}\label{eq:relative-quotient-postnikov-support}
\begin{aligned}
 (q^-)^{-1}\operatorname{PC}^-(q^-\mathcal J^-)
       &=\mathbf D^-_{\mathrm{fg}}(R)_{\operatorname{Supp}\mathcal J^-},\\
 (q^+)^{-1}\operatorname{PC}^+(q^+\mathcal J^+)
       &=\mathbf D^+_{\mathrm{fg}}(R)_{\operatorname{Supp}\mathcal J^+}.
\end{aligned}
\end{equation}
The positive formula also holds internally for finite terms.

If every point of $U$ is closed,
Theorem~\ref{thm:central-relative-scheme} gives
$\mathcal D^-=\mathbf D^-_{\mathrm{fg}}(R)_{W'}$.
For the idempotent $e$ of the complementary factor,
\[
 \operatorname{Supp}(eR)=W',\qquad
 Y\in\mathbf D^+_{\mathrm{fg}}(R)_{W'}
 \ \Longrightarrow\ Y\simeq eR\odot Y\in\mathcal D^+.
\]
The same holds internally. Apply the absolute fourteen conditions
to $R_U=(1-e)R$.

Otherwise choose a homomorphism $R\to V$ over a nonclosed $y\in U$
as in the hypotheses, and
$r\in R$ with image $\varepsilon\pi^c$, $c>0$. Write
$Z=W'\cup V(r)$ and $F=V\otimes_R^{\mathbf L}(-)$. In the negative
category let
\[
 \mathcal J^-=\sqrt{\operatorname{Thick}^{\otimes}
                    (\mathcal D^-\cup\{\operatorname{Cone}(r)\})}.
\]
The functor $F$ sends these generators to uniformly $\pi$-torsion
objects. This class is radical in $\mathbf D^-_{\mathrm{fg}}(V)$,
so the object $K_f(r)$ with unbounded $f$ gives
\[
 K_f(r)\notin\mathcal J^-\subsetneq
 \operatorname{PC}^-(\mathcal J^-)=\mathbf D^-_{\mathrm{fg}}(R)_Z.
\]
The bounded-part theorem and
\eqref{eq:relative-quotient-postnikov-support} disprove (3), (5),
(7). Corollary~\ref{thm:central-relative-growth} disproves (2).

For either positive realization, write
\[
 \mathcal J^+=\operatorname{Thick}^{\operatorname{Hom}}
     (\mathcal D^+\cup\mathbf D^{\mathrm{perf}}(R)_Z),\qquad
 Y_r=\bigoplus_{n\ge0}\operatorname{Cone}(r^{n+1})^\vee[-n].
\]
Here and below the finite-term version is read internally.
For $P\in\mathbf D^{\mathrm{perf}}(R)_Z$, since $y\notin Z$,
\[
 a\in\operatorname{Ann}(\operatorname{id}_P)\setminus y
 \ \Longrightarrow\
 a\operatorname{id}_{B\odot P}=0\quad
       (B\in\mathbf D^-_{\mathrm{fg}}(R)),\qquad u_y(a)\ne0.
\]
Together with \eqref{eq:relative-positive-trait-condition} and
Lemma~\ref{lem:hom-generation-normal-form}, this shows that every
object of $F(\mathcal J^+)$ has a uniform $\pi$-annihilator.
The complex $Y_r$ has finite free terms and is a direct sum of
perfect complexes, whereas
\[
 H^{n+1}F(Y_r)\cong V/\pi^{c(n+1)}V\quad(n\ge0).
\]
Consequently $Y_r\notin\mathcal J^+$, while
\eqref{eq:relative-quotient-postnikov-support} gives
\[
 (q^+)^{-1}\operatorname{PC}^+(q^+\mathcal J^+)
   =\mathbf D^+_{\mathrm{fg}}(R)_Z
   =(q^+)^{-1}\operatorname{PC}^+
                (q^+\mathbf D^+_{\mathrm{fg}}(R)_Z).
\]
This disproves (4), (6), (8) and, internally, (11), (12), (13).
Finally, the negative $K_f(r)$ and positive $Y_r$ become zero
geometrically after inverting $r$, but lie outside
$\operatorname{Thick}(\mathcal D^\alpha\cup\mathcal S^\alpha_r)$:
every object of that thick subcategory has uniformly
$\pi$-annihilated pullback. They disprove (9), (10), (14).
\end{proof}

For any commutative noetherian $R$ and $s\in R$, the family
$\mathcal E=\{\operatorname{Cone}(s)\otimes_R^{\mathbf L}A:
A\in\mathbf D^-_{\mathrm{fg}}(R)\}$ gives
$\mathcal C^\alpha=\mathcal S^\alpha_s$ in all three realizations
and satisfies \eqref{eq:relative-positive-trait-condition}.
Taking $\mathcal E=\varnothing$ and $\Phi=\{1\}$ recovers
Corollary~\ref{thm:artinian-equivalences}.

\end{document}